\let\orilabel\label
\documentclass[onefignum,onetabnum]{siamart190516}
\pdfoutput=1 

\let\label\orilabel

\usepackage{amsfonts}
\usepackage{amsmath}
\usepackage{lmodern}
\usepackage{mathtools}
\usepackage[T1]{fontenc}
\usepackage{graphicx}
\usepackage{epstopdf}
\usepackage{algorithmic}
\usepackage{booktabs}
\usepackage[T1]{fontenc}
\usepackage{tikz,caption}
\usepackage{subcaption}
\usepackage{pgfplots}
\pgfplotsset{compat=newest}
\pgfplotsset{plot coordinates/math parser=false}
\usepackage{url}
\usepackage{color}
\usetikzlibrary{plotmarks} 
\usetikzlibrary{external}
\usepgfplotslibrary{external} 
\tikzset{external/mode=graphics if exists}
\usepackage{hyperref}
\usepackage{cleveref}
\usepackage{multirow}
\usepackage{comment}
\newlength\figureheight
\newlength\figurewidth

\usepackage{xcolor}

\newcommand{\R}{\mathbb{R}}

\newcommand{\hk}{^{\left( d \right)}}
\newcommand{\D}{\,\mathrm{d}}

\def\gmres{{{\sc Gmres}}}

\def\amen{{{\sc AMEn}}}
\def\ttgmres{{{\sc TT-GMRES}}}

\numberwithin{theorem}{section}
\newcommand{\TheTitle}{IETI-based Low-Rank method for PDE-constrained optimization}

\title{{\TheTitle}}

\author{Tom-Christian Riemer\thanks{Technische Universit\"at Chemnitz, Department of Mathematics, Chair of Scientific Computing, 09107 Chemnitz, Germany, \email{tom-christian.riemer@mathematik.tu-chemnitz.de}}
\and
Alexandra B\"unger\thanks{University of British Columbia, Computer Science, Vancouver, BC Canada V6T 1Z4, 604 822 3061, \email{alexandra.buenger@mathematik.tu-chemnitz.de}}
\and 
Martin Stoll\thanks{Technische Universit\"at Chemnitz, Department of Mathematics, Chair of Scientific Computing, 09107 Chemnitz, Germany, \email{martin.stoll@mathematik.tu-chemnitz.de}}
}
\usepackage{amsopn}

\ifpdf
\fi

\definecolor{grey}{rgb}{0.5,0.5,0.5}
\definecolor{darkgreen}{rgb}{0,0.55,0}

\begin{document}

\maketitle

\begin{abstract}
Isogeometric Analysis (IgA) is a versatile method for the discretization of partial differential equations on complex domains, which arise in various applications of science and engineering. Some complex geometries can be better described as a computational domain by a multi-patch approach, where each patch is determined by a tensor product Non-Uniform Rational Basis Splines (NURBS) parameterization. This allows on the one hand to consider the problem of the complex assembly of mass or stiffness matrices (or tensors) over the whole geometry locally on the individual smaller patches, and on the other hand it is possible to perform local mesh refinements independently on each patch, allowing efficient local refinement in regions of high activity where higher accuracy is required, while coarser meshes can be used elsewhere. Furthermore, the information about differing material models or properties that are to apply in a subdomain of the geometry can be included in the patch in which this subdomain is located. For this it must be ensured that the approximate solution is continuous over the entire computational domain and therefore at the interfaces of two (or more) patches. The most promising approach for this problem, which transfers the idea of Finite Element Tearing and Interconnecting (FETI) methods into the isogeometric setup, was the IsogEometric Tearing and Interconnecting (IETI) method, where by introducing a constraints matrix and associated Lagrange multipliers and formulating it into a dual problem, depending only on the Lagrange multipliers, continuity at the interfaces was ensured in solving the resulting system. In this paper we illustrate that low-rank methods based on the tensor-train format can be generalised for a multi-patch IgA setup, which follows the IETI idea.
\end{abstract}

\begin{keywords}
isogeometric analysis, multi-patch, IETI, optimal control, low-rank decompositions, tensor-train format
\end{keywords}

\begin{AMS}
65F10, 65F50, 15A69, 93C20 
\end{AMS}

\section{Motivation}

Isogeometric Analysis (IgA) is a discretization technique used for approximating solutions to a partial differential equation (PDE) defined on a given domain $\Omega$. It was introduced by Hughes, Cottrell and Bazilevs in 2005~\cite{CAD}. In Isogeometric Analysis the problem domain $\Omega$ and the solution space for solving the PDE using a Galerkin approach~\cite{strang} are parameterized by the same spline functions, typically B-splines or NURBS (Non-Uniform Rational Basis Splines). These basis functions are globally defined and have overlapping supports depending on their degrees. As such these discretizations have a higher computational complexity, increasing exponentially with respect to the dimension of the problem~\cite{Mantzaflaris_space_time}, but also allow the relatively easy approximation of domains rather difficult to treat with tradtional finite element methods. One of the major research interests in IgA is to find strategies to overcome the complexity drawback and efficiently assemble the system matrices~\cite{Antolin, Hughes, pan2020fast, pan2021efficient}.

We here follow the idea of Mantzaflaris et al.~\cite{angelos1,angelos2} of using a low-rank tensor method, which exploits the tensor structure of the basis functions and separates the variables of the integrals. As a result the system matrices are then approximated to high accuracy by a sum of Kronecker products of smaller matrices, which are assembled via univariate integration. We here rely on the method of~\cite{BuengerDolgovStoll:2020} where the assembly is carried out using an interpolation step and a low-rank representation of the resulting coefficient tensor. The authors there combine the low-rank method of Mantzaflaris et al.\ with low-rank tensor-train (TT) calculations~\cite{osel-tt-2011, DoOs-dmrg-solve-2011}. Exploiting the tensor product nature of the arising interpolation, we can calculate a low-rank TT approximation without prior assembly of the full coefficient tensor by means of the Alternating Minimal Energy (AMEn) method~\cite{amen}. Our goal for this paper is to extend this technique to the case of a multi-patch domain discretized using IgA. This poses the problem that the approximations generated for the corresponding problems can show discontinuities at the interfaces of these patches. We overcome this by transferring the idea of the IETI method from~\cite{KleissPechsteinJuettlerTomar:2012} to the low-rank tensor setup. We also want to test our technique on the following problems. We consider the low-rank solution of the elliptic problem defined by Poisson's equation equipped with homogeneous Dirichlet boundary conditions
\begin{equation}
\label{eq:forward}
\begin{aligned}
-\Delta y &= f \quad \mbox{ in } \Omega, \\
y & = 0  \quad \mbox{ on } \partial \Omega.
\end{aligned}
\end{equation}
Here $f$ is some source function and $\Omega$ is a given mulit-patch geometry parameterized by B-splines or NURBS. The second problem we consider is an optimization problem where the heat equation becomes the constraint of an objective function that we want to minimize, i.e.,
\begin{align}
\min_{y,u} \quad \frac{1}{2}\int_0^T \! \int_\Omega (y - & \hat{y})^2  \D x\D t + \frac{\alpha}{2}\int_0^T \, \int_\Omega u^2 \D x \D t \label{eq:optimization1} \\
\mbox{s.t.}\quad y_t - \Delta y & = u \quad \mbox{ in } (0,T)\times \Omega, \label{eq:optimization2}\\
y & = 0 \quad \mbox{ on } (0,T)\times \partial \Omega, \label{eq:optimization3}\\
y & = y_0 \quad \mbox{ on } \Omega \mbox{ for } t=0, \label{eq:optimization4}
\end{align}
with a desired state $\hat{y}$ and control $u$ on a mulit-patch geometry $\Omega$. The discretization of~\eqref{eq:optimization1} to~\eqref{eq:optimization4} in this paper will be performed by IgA. Tensor techniques for IgA have shown promising results in many areas and we refer to~\cite{montardini2023lowsp} for the single-patch case and for the multi-patch case to~\cite{montardini2024lowrank}. In the latter the authors use Tucker tensors for the low-rank approximation and focus on the forward elasticity simulation on conforming (or fully matching) multi-patch geometries. In this paper we focus on the approximation via the tensor-train format and provide approaches for the solution of the optimal control problem, also on nonconforming geometries. Nevertheless, their method and ours are similar in spirit by aiming at breaking the curse of dimensionality by relying on low-rank tensor formats.

The paper is structured as follows: In the preliminaries we first discuss low-rank tensor formats in subsection~\ref{subsection:Low-rank_tensor_format}, in subsection~\ref{subsection:Low-rank_IgA} we introduce the basics of IgA and discuss how the tensor-train format can be used so that the system matrices or tensors can be assembled low-rank, in subsection~\ref{subsection:Multi-patch_IgA} we present our multi-patch IgA notation and state the general problem. In section~\ref{section:jump_tensors}, we present how the idea of the IETI method~\cite{KleissPechsteinJuettlerTomar:2012} can be generalised for the tensor setup and how the so-called jump tensors can be defined. Then we explain in section~\ref{section:IETI_method} how the resulting low-rank IETI method works to generate continuous low-rank approximations over multi-patch geometries. In section~\ref{section:optimization}, we show how this method can be used to find an approximation of the optimization problem described by~\eqref{eq:optimization1} to~\eqref{eq:optimization4}. The results of our numerical experiments are presented in section~\ref{section:numerics}. In section~\ref{section:conclusion}, we summarise the insights and results, concluding our work.

\newpage

\section{Preliminaries} \label{section:basics}
\subsection{Low-rank tensor format}
\label{subsection:Low-rank_tensor_format}
The most well-known technique for low-rank approximations is the singular value decomposition, illustrated for a matrix $W \in \R^{n_1\times n_2}$ as
\begin{equation}
\label{eq:SVD}
W = U \Sigma V^\top \approx \sum_{r=1}^R u_r \sigma_r v_r^\top = \sum_{r=1}^R (u_r \sqrt{\sigma_r}) \otimes (v_r \sqrt{\sigma_r}).
\end{equation}
with $U\in \R^{n_1\times n_1}$, $V \in \R^{n_2\times n_2}$ with their columns denoted by $u_r$ and $v_r$, and $\Sigma \in \R^{n_1 \times n_2}$ is the rectangular matrix holding the sorted singular values $\sigma_i$, $i=1,\ldots,\min(n_1,n_2)$ on its main diagonal. The best low-rank approximation is obtained by the truncated SVD where we truncate all singular values below some given threshhold resulting in a rank-$R$ approximation, where $R$ is the number of used singular values and therefore the number of summands in~\eqref{eq:SVD}.

In the high-dimensional case we need low-rank tensor approximation of a $D$-dimensional tensor. Such approximations are given by, e.g., the higher-order singular value decomposition (HOSVD)~\cite{mlsvd}, or a canonical polyadic decomposition (CP)~\cite{CPD}. However, the approximation problem in the CP format is typically ill-posed~\cite{desilva-2008} and might be numerically unstable. The HOSVD (known also as the Tucker format) still contains the curse of dimensionality as it relies on the dimension of the original tensor. We switch to the more robust tensor-train (TT) decomposition~\cite{osel-tt-2011} in this paper also given the availability of appropriate methods within a robust software framework.

A tensor $W \in \mathbb{R}^{n_1 \times \ldots \times n_D}$ is given in the TT format if it is written as
\begin{align}
\label{eq:tt}
W(i_1,\ldots,i_D) = W^{\left(1\right)}(i_1) \cdots W^{\left( D \right)}(i_D),
\end{align}
where $W^{\left(d\right)}(\cdot) \in \mathbb{R}^{R_{d-1} \times n_d \times R_d}$ are the TT cores, which can be understood as parameter dependent matrices $W^{\left(d\right)}(i_d)$, $i_d = 1, \ldots, n_d$, of size $R_{d-1}\times R_d$ with $R_0 = R_D = 1$~\cite{osel-tt-2011}. The TT format can be rewritten into a canonical representation as
\begin{equation} 
\label{eq:TTranks}
W = \sum_{r_1=1}^{R_1}\cdots \sum_{r_D=1}^{R_D} \bigotimes_{d=1}^D W^{\left(d\right)}(r_{d-1},:,r_d).
\end{equation}

\subsection{Low-rank IgA}
\label{subsection:Low-rank_IgA}
Isogeometric analysis allows to represent a geometry exactly using a set of B-splines or NURBS~\cite{piegl} and by using the same basis functions for the solution space of a PDE on this geometry lies at the heart of the IgA method and its success in scientific computing~\cite{CAD,nguyen2015isogeometric}. We here briefly review some important properties of the method with a focus on deriving the discretized equations.

A set of $n$ B-splines is uniquely defined by its degree $p\in \mathbb{N}_0$ and the knot vector $\xi= \{ \xi_1, \ldots, \xi_{n+p+1} \}$ with
\begin{equation} 
0= \xi_1 = \cdots = \xi_{p+1} < \xi_{p+2} \leq \cdots \leq \xi_n < \xi_{n+1} = \cdots = \xi_{n+p+1} = 1, 
\label{eq:knotVector}
\end{equation}
where the end knots appear $p+1$ times and for all other knots, duplicate appearances are allowed up to multiplicity $p$. Here $n \in \mathbb{N}$ denotes the number of B-splines $\beta_{i,p}$, with $i=1,\ldots,n$. 

For each knot vector $\xi$ as in~\eqref{eq:knotVector}, the according B-splines $\beta_{i,p}$ of degree $p$, with $i = 1,\ldots, n$, are uniquely defined by a recursion formula. The resulting B-splines $\beta_{i,p}$ have the local support $[\xi_i,\xi_{i+p+1}]$. We use $\mathbb{S}_\xi^p$ to denote the spline space spanned by the B-splines with degree $p$ and knot vector $\xi$ and we refer to the basis functions as $\beta_i \in \mathbb{S}_\xi^p$. In order to increase the accuracy of the numerical approximation a refinement strategy based on knot insertion is often applied and we refer to~\cite{CAD,piegl}.  

For $D$-dimensional geometries we use tensor products of univariate spline spaces considering $D$ different univariate spline spaces $\mathbb{S}_{\xi_d}^{p_d}$, where one can assume that each space has the degree $p_d$ and an individual knot vector $\xi_d$, with $d = 1,  \ldots ,D$. Here, $\hat{x}\hk \in [0,1]$ are the $1D$ variables and $\beta_{1}\hk, \ldots, \beta_{n_d}\hk$ the basis functions. The resulting spline space is then denoted by $\mathbb{S}_D =  \mathbb{S}_{\xi_1}^ {p_1} \otimes \ldots \otimes \mathbb{S}_{\xi_D}^{p_D}$. For simplicity we assume further that $p_d = p \; \forall d$ and the index $p_d$ will be omitted from $\mathbb{S}_{\xi_d}^{p}$ in the remainder due to better readability. The basis functions of $\mathbb{S}_D$ are denoted by
\begin{equation*}
\beta_\mathbf{i}  \left( \hat{x} \right)  = \prod_{d=1}^D \beta_{i_d}\hk  \left( \hat{x}\hk \right),   
\end{equation*} 
with multi-index $\mathbf{i} \in \mathbf{I} = \left\{ \left(i_1, \ldots, i_D \right) \, \colon \, i_d \in \left\{ 1,\ldots, n_d \right\}, \, d=1,\ldots, D \right\}$ and variables $\hat{x} = \left[ \hat{x}^{(1)}, \ldots, \hat{x}^{(D)} \right]^\top \in [0,1]^D$.
All multivariate basis functions evaluated at a point $\hat{x} \in [0,1]^D$ can be written as a tensor product
\begin{equation*}
B \left( \hat{x} \right) = \bigotimes_{d=1}^D B^{\left( d \right)} \left( \hat{x}^{\left( d \right)} \right) \in \R^{n_1 \times \ldots \times n_D},
\end{equation*}
where $B^{\left(d\right)} \left( \hat{x}^{\left(d\right)} \right) = \left[ \beta_{1_d}^{\left(d\right)} \left( \hat{x}^{\left(d\right)} \right), \ldots, \beta_{n_d}^{\left( d \right)} \left( \hat{x}^{\left( d \right)} \right) \right]^\top \in \mathbb{R}^{n_d}$ is a vector holding the univariate B-splines in dimension $d = 1, \ldots, D$.

To use these functions for solving a PDE on the domain $\Omega \subset R^D$ we need a B-spline geometry mapping $G \colon \hat{\Omega} \rightarrow \Omega$ from the $D$-dimensional unit cube $\hat{\Omega}:=[0,1]^D$ onto $\Omega$. This is given by 
\begin{equation} \label{eq:geometryMapping}
G \left( \hat{x} \right)  = \sum_{\mathbf{i} \in \mathbf{I}} C_\mathbf{i} \beta_{\mathbf{i}} \left( \hat{x} \right)  = C:B \left( \hat{x} \right) ,
\end{equation}
where $C_{\mathbf{i}} \in \R^D$ are the control points. All control points and evaluations of the B-splines are organised in the tensors $C \in \R^{D \times n_1 \times \ldots \times n_D}$ and $B \left( \hat{x} \right) \in \R^{n_1\times \ldots \times n_D},$ respectively. Here, $:$ denotes the Frobenius product. To overcome some limitations of the B-spline approach NURBS (Non-uniform rational B-splines) have been used~\cite{NURBS} quite extensively but will not be discussed further here.

The discretization of the PDE is usually obtained from a weak formulation where we compute approximations of $y \in H^1_0(\Omega)$ with discrete functions $y_h \in V_h \subset H^1_0(\Omega)$ using B-splines. In IgA, the same splines that are used in the construction of the geometry mapping~\eqref{eq:geometryMapping} are used to parameterize the solution space, i.e. $V_h = \operatorname{span} \left\{ \hat{\beta}_\mathbf{i} : = \beta_\mathbf{i} \circ G^{-1} \, : \, \mathbf{i} \in \mathbf{I}_0 \right\} \subset H^1_0(\Omega)$ with basis functions $\beta_{\mathbf{i}} \in \mathbb{S}_D$ and an index set $\mathbf{I}_0 = \left\{ \left(i_1, \ldots, i_D \right) \, \colon \, i_d \in \left\{ 2,\ldots, n_{d}-1 \right\}, \, d=1,\ldots, D \right\} \subset \mathbf{I}$ in which the first and last index of $\mathbf{I}$ in each dimension are omitted, since the remaining splines with index in $\mathbf{I} \setminus \mathbf{I}_0$ are zero for homogeneous Dirichlet conditions. To improve readability, we make an index shift so that $\mathbf{I}_0 = \left\{ \left(i_1, \ldots, i_D \right) \, \colon \, i_d \in \left\{ 1,\ldots, \tilde{n}_{d} \right\}, \, d=1,\ldots, D \right\}$. The functions $y_h \in V_h$ are linear combinations of the basis functions $y_h = \sum_{\mathbf{i} \in \mathbf{I}_0} y_{\mathbf{i}} ( \beta_{\mathbf{i}} \circ G^{-1})$ with coefficients $y_\mathbf{i} \in \R$. The tensor product structure of $\mathbb{S}_D$ induces a tensor product structure of the solution space $V_h$, since each basis function $\hat{\beta}_\mathbf{i} \in V_h$, $\mathbf{i} = \left(i_1, \ldots, i_D \right) \in \mathbf{I}_0$, can be represented as
\begin{equation}
\label{eq:solution_tensor_product}
\hat{\beta}_\mathbf{i} \left( x \right) = \beta_\mathbf{i} \left({G}^{-1} \left(x\right)\right) = \beta^{\left(1\right)}_{i_1} \left( {{G}^{-1} \left(x\right)}^{\left(1\right)} \right) \cdots \beta^{\left(D\right)}_{i_D} \left( {{G}^{-1} \left(x\right)}^{\left(D\right)} \right),
\end{equation}
where ${{G}^{-1} \left(x\right)}^{\left(d\right)} \in \left[0, 1 \right]$ is the $d$-th component of the inverse of the geometry mapping and the $\beta^{\left(d\right)}_{i_d} \in \mathbb{S}_{\xi_d}$, $d = 1, \ldots, D$, are the univariate splines.\\

The space $V_h$ is now used for the Galerkin discretization, resulting in the discrete mass and stiffness terms
\begin{align*}
a_{m}(u_h, v_h) &= \int_\Omega u_h(x) v_h(x) \D x = \int_{\hat{\Omega}} \sum_{\mathbf{i} \in \mathbf{I}_0} u_\mathbf{i} \beta_\mathbf{i} \left( \hat{x} \right)  \sum_{\mathbf{j} \in \mathbf{I}_0} v_\mathbf{j} \beta_\mathbf{j} \left( \hat{x} \right)  \omega \left( \hat{x} \right)  \D \hat{x}, \\
a_{s}(u_h,v_h) &= \int_{\Omega} \nabla u_h(x) \cdot \nabla v_h(x) \D x = \int_{\hat{\Omega}} \left(Q \left( \hat{x} \right)  \sum_{\mathbf{i} \in \mathbf{I}_0} u_\mathbf{i} \nabla \beta_\mathbf{i} \left( \hat{x} \right) \right) \cdot\sum_{\mathbf{j} \in \mathbf{I}_0} v_\mathbf{j} \nabla \beta_\mathbf{j} \left( \hat{x} \right)  \D \hat{x},
\end{align*}
for $u_h, v_h \in V_h$, with the additional terms stemming from the domain transformation,
\begin{align*}
\omega \left( \hat{x} \right)  &= \lvert \det \nabla G \left( \hat{x} \right) \rvert && \in \R, \\
Q \left( \hat{x} \right)  &= \left( \nabla G \left( \hat{x} \right) ^{T} \nabla G \left( \hat{x} \right) \right)^{-1} \lvert \det \nabla G \left( \hat{x} \right) \rvert && \in \R^{D\times D},
\end{align*}
as introduced in~\cite{angelos1}. The corresponding mass and stiffness terms can be written in tensor form, i.e., the mass tensor
\begin{equation} 
\label{eq:mass_tensor}
M = \int_{\hat{\Omega}} \omega \left( \hat{x} \right) \, B \left( \hat{x} \right) \otimes B \left( \hat{x} \right) \D \hat{x} \quad \in \R^{\left(\tilde{n}_1, \ldots, \tilde{n}_D \right) \times \left( \tilde{n}_1 , \ldots, \tilde{n}_D \right)}.
\end{equation} 
Similarly, we can write the stiffness tensor as
\begin{equation}
\label{eq:stiffness_tensor}
K = \int_{\hat{\Omega}} \left[ Q \left( \hat{x} \right) \cdot \left( \nabla \otimes B \left( \hat{x} \right) \right) \right] \cdot \left( \nabla \otimes B \left( \hat{x} \right) \right)  \D \hat{x} \quad \in \R^{\left(\tilde{n}_1, \ldots, \tilde{n}_D \right) \times \left( \tilde{n}_1 , \ldots, \tilde{n}_D \right)}.
\end{equation} 

The computation and storage of~\eqref{eq:mass_tensor} and~\eqref{eq:stiffness_tensor} can be extremely expensive due to multi-dimensional quadrature and the overlapping support of B-splines with high degrees. But it has been observed that these tensors can be well approximated in low-rank tensor formats~\cite{angelos1,BuengerDolgovStoll:2020} based on a low-rank approximation of the coupling terms in the integral. For that we approximate the arising multi-dimensional integrals as products of univariate integrals. The ingredients for the mass and stiffness tensors are all univariately defined, except for the weight functions $\omega \left( \hat{x} \right)$ and $Q \left( \hat{x} \right)$, which are determined by the geometry mapping~\eqref{eq:geometryMapping} but are not separable into one-dimensional factors. We therefore interpolate these weight functions by a combination of univariate B-splines of higher order $\hat{B} \left(\hat{x}\right) \in \mathbb{R}^{\left( \hat{n}_1, \ldots, \hat{n}_D \right)}$, i.e.
\begin{equation}
\label{eq:weightfunction_linear_system}
\omega \left( \hat{x} \right) \approx W \colon \hat{B} \left(\hat{x}\right).
\end{equation}
For that we follow the approach introduced in~\cite{BuengerDolgovStoll:2020} by computing a low-rank approximation in TT format~\eqref{eq:TTranks} of the coefficient tensor $W \in \mathbb{R}^{\left( \hat{n}_1, \ldots, \hat{n}_D \right)}$ 
\begin{equation*}
W_R : = \sum_{r_1=1}^{R_1}\cdots \sum_{r_D=1}^{R_D} \bigotimes_{d=1}^D W^{\left(d\right)}_R (r_{d-1},:,r_d) = \sum_{r = 1}^{R} \bigotimes^{D}_{d=1}  w_r^{(d)} \approx W,
\end{equation*}
where $w_r^{(d)} \in \mathbb{R}^{\hat{n}_d}$ and $R = R_1 \cdots R_D$. With this we get a low-rank representation of the weight function,
\begin{equation*}
\omega \left( \hat{x} \right) \approx W_R \colon \hat{B} \left(\hat{x}\right) = \sum_{r=1}^R \prod_{d=1}^D w_r^{(d)}\cdot \hat{B}^{(d)} \left( \hat{x}^{(d)} \right),
\end{equation*}
where $\hat{B}^{(d)} \left( \hat{x}^{(d)} \right) \in \mathbb{R}^{\hat{n}_d}$ denotes the vector holding all univariate basis functions evaluated in $\hat{x}^{(d)} \in \left[0, 1\right]$. As a result the integrands are separable and we can write the mass tensor~\eqref{eq:mass_tensor} as a sum of tensor products of small univariate mass matrices
\begin{equation} 
\label{eq:massFinal}
\begin{aligned}
M &= \sum_{r=1}^R \bigotimes_{d=1}^D \int_{0}^{1} \left( w_r^{(d)}\cdot \hat{B}^{(d)} \left( \hat{x}^{(d)} \right) \right) \, B^{(d)} \left( \hat{x}^{(d)} \right) \otimes B^{(d)} \left( \hat{x}^{(d)} \right) \D \hat{x}^{(d)}  \\
&= \sum_{r=1}^R \bigotimes_{d=1}^D M_r^{\left( d \right)}.
\end{aligned}
\end{equation}
The same procedure can be applied to each entry of $Q \left( \hat{x} \right)$ such that we get a low-rank tensor representation of~\eqref{eq:stiffness_tensor} as
\begin{equation} 
\label{eq:stiffnessFinal}
K = \sum_{k,l=1}^D \sum_{r=1}^R \bigotimes_{d=1}^D K_{k,l,r}^{\left( d \right)}.
\end{equation}
We refer to~\cite{BuengerDolgovStoll:2020} for details and to the codes on our website \cite{BuengerCode}.

\subsection{Multi-patch IgA} 
\label{subsection:Multi-patch_IgA}

In a multi-patch setting we assume that the geometric shape $\Omega \subset \R^D$ can be decomposed into $N_P$ many single-patch NURBS parameterizations, such as~\eqref{eq:geometryMapping}, i.e.\ 
\begin{equation*}
G^{\left( j \right)} \left( \hat{\Omega} \right) = \Omega^{\left( j \right)} \subset \Omega, \quad j = 1, \ldots, N_P,
\end{equation*}
such that
\begin{equation*}
\overline{\Omega} = \bigcup^{N_P}_{j = 1} \overline{\Omega^{\left( j \right)}} \quad \text{and} \quad \Omega^{\left( j \right)} \cap \Omega^{(k)} = \emptyset.
\end{equation*}
We note that the parameter space $\hat{\Omega}:=[0,1]^D$ is the same for each parameterization $G^{\left( j \right)}$. The multi-patch geometries considered in this paper are all $3$-dimensional, so from now on we set $D=3$ to simplify the notation. Further we assign all variables belonging to the patch $\Omega^{\left( j \right)}$, such as basis functions, control points, index sets, etc.\ a superscript $\left( j \right)$. 

As in the single-patch case, the same splines to represent the geometry $\Omega$ are used to approximate the solution of the underlying PDE problem. Since each patch $\Omega^{\left( j \right)}$ has its own tensor product spline space $\mathbb{S}^{\left( j \right)}_3 = \mathbb{S}^{\left( j \right)}_{\xi_1} \otimes \mathbb{S}^{\left( j \right)}_{\xi_2} \otimes \mathbb{S}^{\left( j \right)}_{\xi_3}$ with parameterization $G^{\left( j \right)} \colon \hat{\Omega} \to \Omega^{\left( j \right)}$, we define a local solution space for each patch via
\begin{equation*}
V^{\left( j \right)}_h = \operatorname{span} \left\{ \hat{\beta}^{\left( j \right)}_\mathbf{i} := \beta^{\left( j \right)}_\mathbf{i} \circ {G^{\left( j \right)}}^{-1} \, : \, \mathbf{i} \in \mathbf{I}^{\left( j \right)}_0 \right\} \subset H^1(\Omega^{\left( j \right)}),
\end{equation*}
where $\mathbf{I}^{\left( j \right)}_0 \subset \mathbf{I}^{\left( j \right)}$ contains the indices of all splines of patch $\Omega^{\left( j \right)}$ whose support does not lie on $\partial \Omega \cap \partial \Omega^{\left( j \right)}$. As in the single-patch case, the splines of the remaining indices can be considered as zero. We assume that this set of indices has a tensor structure, i.e. $\mathbf{I}^{\left( j \right)}_0 = \left\{ \left(i^{\left( j \right)}_1, i^{\left( j \right)}_2, i^{\left( j \right)}_3 \right) \, \colon \, i^{\left( j \right)}_d \in \left\{ 1, \ldots, \tilde{n}^{\left( j \right)}_{d} \right\}, \, d=1,2,3 \right\} \subset \mathbb{N}^{\tilde{n}^{\left( j \right)}_1 \times \tilde{n}^{\left( j \right)}_2 \times \tilde{n}^{\left( j \right)}_3}$. We note that the solution space $V^{\left( j \right)}_h$, $j = 1, \ldots, N_p$, has a tensor product structure induced by $\mathbb{S}^{\left( j \right)}_3$, since all basis functions can be written as in~\eqref{eq:solution_tensor_product}. 

The space of functions on $\Omega$ which are locally in $V^{\left( j \right)}_h$ is denoted by 
\begin{equation*}
\Pi V_h = \left\{ y \in \mathcal{L}^2 \left( \Omega \right) \colon \; y|_{\Omega^{\left( j \right)}} \in V^{\left( j \right)}_h,  \; \forall j = 1, \ldots, N_{p} \right\}.
\end{equation*}
Each function $y_h \in \Pi V_h$, $y_h \colon \Omega \to \mathbb{R}$, can be represented patch-wise by a linear combination of the basis functions of the corresponding patch $V^{\left( j \right)}_h$ with coefficients $y^{\left( j \right)}_\mathbf{i} \in \R$, which from now on are referred to as degrees of freedom (DoFs), i.e.
\begin{equation}
\label{eq:discrete_function}
y_h|_{\Omega^{\left( j \right)}} \left( \tilde{x} \right) = \sum_{\mathbf{i} \in \mathbf{I}^{\left( j \right)}_0} y^{\left( j \right)}_{\mathbf{i}} \, \hat{\beta}^{\left( j \right)}_{\mathbf{i}} \left( \tilde{x} \right), \quad \forall j = 1, \ldots, N_P.
\end{equation}
The set $\mathbf{I}^{\left( j \right)}_0$ can therefore be seen as the set of DoFs for each patch $\Omega^{\left( j \right)}$, $j = 1, \ldots, N_P$. \\

The multi-patch approach is used when the geometric domain $\Omega$ cannot be parameterized by a single geometry mapping~\eqref{eq:geometryMapping} but it also enables us to assume different material models and element types on different patches (cf.~\cite{CAD}) or to undertake a patch-wise local refinement by using different rich bases for the solution spaces of different patches. We note that the patches coincide with the non-overlapping subdomains of a FETI-like method~\cite{KleissPechsteinJuettlerTomar:2012}. The challenge in computing an approximation for a PDE problem on a multi-patch geometry $\Omega$ lies in the fact that discrete functions $y_h \colon \Omega \to \mathbb{R}$ defined by~\eqref{eq:discrete_function} are, in general, discontinuous across the patch interfaces. In the following we denote the interface of two patches $\Omega^{\left( j \right)}$ and $\Omega^{\left( k \right)}$ by 
\begin{equation*}
\Gamma^{\left(j, k\right)} = \partial \Omega^{\left( j \right)} \cap \partial \Omega^{\left( k \right)}. 
\end{equation*}
The set of the index-tupels of all interfaces that are not empty is denoted by 
\begin{equation*}
\mathcal{C} = \left\{ (j, k) \in \left\{ 1, \ldots, N_P \right\}^2 \, \colon  \, \Gamma^{\left(j, k\right)} \neq \emptyset, \, j < k \right\}.
\end{equation*}
The condition $j<k$ ensures that each interface is only counted once in $\mathcal{C}$. 

In this paper we assume that the computational domain $\Omega$ is represented as a collection of several patches connected along their interfaces with $C^0$-continuity. Since the parameter space is the unit cube $\hat{\Omega}=[0,1]^3$ for each patch, we further assume that each interface $\Gamma^{\left(j, k\right)}$ of two patches $\Omega^{\left( j \right)}$ and $\Omega^{\left( k \right)}$ is always a $2$-dimensional surface, which is the image of one entire side of the six sides of the unit cube $\hat{\Omega}$ under both parameterizations $G^{\left( j \right)}$ and $G^{\left( k \right)}$. For the sake of simplicity, we assume that the two patches which are connected via an interface $\Gamma^{\left(j, k\right)}$ have the same orientation in the Cartesian coordinate system, i.e.\ if we number the sides of the unit cube $\hat{\Omega}$ like a dice for each patch (e.g.\ the side $\{ 0 \} \times \left[0, 1 \right] \times \left[0, 1 \right]$ is referred to as side $1$ and $\{ 1 \} \times \left[0, 1 \right] \times \left[0, 1 \right]$ as side $6$) and for patch $\Omega^{\left( j \right)}$ the interface $\Gamma^{\left(j, k\right)}$ is the image under $G^{\left( j \right)}$ of side $1$, then for patch $\Omega^{\left( k \right)}$ the interface $\Gamma^{\left(j, k\right)}$ is the image under $G^{\left( k\right)}$ of side $6$. This means that two opposite sides of the dice are always the sides of the patches that form their interface. We also assume that there can only be one interface between two patches. 

We say that for $(j, k) \in \mathcal{C}$ the two patches $\Omega^{\left( j \right)}$ and $\Omega^{\left( k \right)}$ are connected in dimension $d_1 \in \{ 1, 2, 3 \}$ by the interface $\Gamma^{\left(j, k\right)}$, if both parameterizations $G^{\left( j \right)}$ and $G^{\left( k \right)}$ are not fixed in the parameters of the dimensions $d_2, d_3 \in \{ 1, 2, 3 \} \setminus \{ d_1 \}$, $d_2 \neq d_3$, when mapping to the interface $\Gamma^{\left(j, k\right)}$, which means that the dimensions $d_2$ and $d_3$ span the interface $\Gamma^{\left(j, k\right)}$. In the following, triplets of the variables that depend on these dimensions are written in the order specified by the additional index, so that it can be understood as a correctly permuted variant of the triplet, i.e.\ let $d_1 = 3$, $d_2 = 1$, $d_3 = 2$, then $\left( \tilde{n}^{\left(j\right)}_{d_1}, \tilde{n}^{\left(j\right)}_{d_2}, \tilde{n}^{\left(j\right)}_{d_3} \right)$ is equivalent to $\left( \tilde{n}^{\left(j\right)}_{1}, \tilde{n}^{\left(j\right)}_{2}, \tilde{n}^{\left(j\right)}_{3} \right)$. \\
In the following for $(j, k) \in \mathcal{C}$ we will use 
\begin{equation*}
\mathbf{I}_{\Gamma} (j, k) = \left\{ \mathbf{i} \in \mathbf{I}_{0}^{\left( j \right)} \, \colon \, \operatorname{supp} \left( \hat{\beta}^{\left( j \right)}_{\mathbf{i}} \right) \cap \Gamma^{\left(j, k\right)} \neq \emptyset \right\}
\end{equation*}
to denote the set of indices of basis functions on $\Omega^{\left( j \right)}$, whose support intersects with the interface $\Gamma^{\left(j, k\right)}$ and for $\mathbf{i} \in \mathbf{I}_{\Gamma} (j, k)$ we say that $y^{\left( j \right)}_{\mathbf{i}}$ is associated with the interface $\Gamma^{\left(j, k\right)}$. The definition of $\mathbf{I}_{\Gamma}(k, j)$ is analogous for $(j, k) \in \mathcal{C}$, i.e.\ for $\mathbf{m} \in \mathbf{I}_{\Gamma} (k, j)$, $y^{\left( k \right)}_{\mathbf{m}}$ is associated with $\Gamma^{\left(j, k\right)}$. Because of our assumption about the orientation of the patches and the tensor product structure of the splines, these DoFs can be easily identified. For example, let for $(j, k) \in \mathcal{C}$ the patches $\Omega^{\left( j \right)}$ and $\Omega^{\left( k \right)}$ be connected in dimension $d = 2$, such that the interface $\Gamma^{\left(j, k\right)}$ is located on side $2$ of patch $\Omega^{\left( j \right)}$ and correspondingly on side $5$ of patch $\Omega^{\left( k \right)}$, then 
\begin{equation*}
\begin{aligned}
\mathbf{I}_{\Gamma} (j, k) & = \left\{ \left( i^{\left( j \right)}_1, 1, i^{\left( j \right)}_3 \right) \, \colon \, i^{\left( j \right)}_d = 1, \ldots, \tilde{n}^{\left( j \right)}_d, \, d \in \left\{1, 3 \right\} \right\} \subset \mathbb{N}^{\tilde{n}^{\left( j \right)}_1 \times 1 \times \tilde{n}^{\left( j \right)}_3}, \\
\mathbf{I}_{\Gamma} (k, j) & = \left\{ \left( i^{\left( k \right)}_1, \tilde{n}^{\left( k \right)}_2, i^{\left( k \right)}_3 \right) \, \colon \, i^{\left( k \right)}_d = 1, \ldots, \tilde{n}^{\left( k \right)}_d, \, d \in \left\{1, 3 \right\} \right\} \subset \mathbb{N}^{\tilde{n}^{\left( k \right)}_1 \times 1 \times \tilde{n}^{\left( k \right)}_3}.
\end{aligned}
\end{equation*}
We see, if the two patches $\Omega^{\left( j \right)}$ and $\Omega^{\left( k \right)}$ are connected in dimension $d_1 \in \{1, 2, 3 \}$ by the interface $\Gamma^{\left(j, k\right)}$, then $\mathbf{I}_{\Gamma} (j, k)$ can be understood as the set of the $\tilde{n}^{\left( j \right)}_{d_2} \tilde{n}^{\left( j \right)}_{d_3}$ many DoFs $y^{\left(j\right)}_{\mathbf{I}_{\Gamma} (j, k)} \in \mathbb{R}^{1 \times \tilde{n}^{\left( j \right)}_{d_2} \times \tilde{n}^{\left( j \right)}_{d_3}}$ of patch $\Omega^{\left( j \right)}$ that lie on the interface $\Gamma^{\left(j, k\right)}$.

\section{Jump tensors}
\label{section:jump_tensors}
The strategy of the IETI method~\cite{KleissPechsteinJuettlerTomar:2012} for finding an approximation for~\eqref{eq:forward} is to determine the coefficients $\mathbf{y}$ of the discrete approximation $y_h \in \Pi V_h$ by solving the following saddle point formulation
\begin{equation}
\begin{bmatrix}
\mathbf{K} & \mathbf{A}^\top \\
\mathbf{A} & 0
\end{bmatrix}
\begin{bmatrix}
\mathbf{y} \\
\boldsymbol{\lambda}
\end{bmatrix}
= 
\begin{bmatrix}
\mathbf{f} \\
0
\end{bmatrix},
\label{eq:SaddlePointFormulation}
\end{equation}
where $\mathbf{K}$ is a block diagonal matrix having the local stiffness matrices of each patch on its diagonal, $\mathbf{A}$ is a so-called jump matrix through which the $C^0$-continuity is enforced by linear constraints, $\mathbf{f}$ is the source vector and $\boldsymbol{\lambda}$ the corresponding Lagrange multipliers. To transfer this idea to the tensor setup, we think of the system~\eqref{eq:SaddlePointFormulation} as a block system of tensors, where each block has again a block structure. This means 
\begin{equation}
\mathbf{K} = 
\begin{bmatrix}
K^{\left( 1 \right)} & & \\
& \ddots & \\
& & K^{\left( N_P \right)}
\end{bmatrix}
\label{eq:StiffnessMatrixMultiPatch}
\end{equation}
is a block diagonal tensor and its diagonal blocks $K^{\left( j \right)}$, $j = 1, \ldots, N_P$, are $3$-\\
dimensional stiffness tensors defined by~\eqref{eq:stiffnessFinal}, corresponding to the bilinear form on each patch $\Omega^{\left( j \right)}$. The tensor $\mathbf{y} = \left[ y^{\left( 1 \right)}, \ldots, y^{\left( N_P \right)} \right]^\top$ is the unique representation of $y_h \colon \Omega \to \mathbb{R}$, whose blocks $y^{\left( j \right)} \in \mathbb{R}^{\tilde{n}^{(j)}_1 \times \tilde{n}^{(j)}_2 \times \tilde{n}^{(j)}_3}$ are the local tensors with the real-valued coefficients of each patch $\Omega^{\left( j \right)}$ for~\eqref{eq:discrete_function}. The source tensor $\mathbf{f}$ has the same structure as $\mathbf{y}$.

The continuity of the approximation $y_h$ is ensured in~\eqref{eq:SaddlePointFormulation} by a so-called \textit{jump tensor} $\mathbf{A}$ and the Lagrange multipliers $\boldsymbol{\lambda} = \left[ \lambda^{\left( 1 \right)}, \ldots, \lambda^{\left( \lvert \mathcal{C} \rvert \right)} \right]^\top$. The jump tensor $\mathbf{A}$ is also in block structure and has horizontally $N_P$ many block columns and vertically $\lvert \mathcal{C} \rvert$ many block rows, one for each interface. In each block row of $\mathbf{A}$, all blocks are zero tensors except for two blocks. Let $\left(j, k\right) \in \mathcal{C}$, then by applying $\mathbf{A}$ to $\mathbf{y}$, the $C^0$-continuity of $y_h \colon \Omega \to \mathbb{R}$ on the interface $\Gamma^{\left(j, k\right)}$ is enforced on the DoFs $y^{\left( j \right)}$ and $y^{\left( k \right)}$ through the corresponding block row of $\mathbf{A}$ with the two non-zero tensors $A^{\left( j \right)}_{\left(j, k\right)}$ and $A^{\left( k \right)}_{\left(j, k\right)}$ via 
\begin{equation}
\label{eq:DoFs_subtraction}
A^{\left( j \right)}_{\left(j, k\right)} \cdot y^{\left( j \right)} - A^{\left( k \right)}_{\left(j, k\right)} \cdot y^{\left( k \right)} = 0,
\end{equation}
where $\cdot$ denotes the contracted product over the dimensions $\tilde{n}^{(j)}_1 \times \tilde{n}^{(j)}_2 \times \tilde{n}^{(j)}_3$. The idea is that the DoFs of one patch should be expressed as a linear combination of the DoFs of the other patch on the interface $\Gamma^{\left(j, k\right)}$.
We note that other boundary conditions are not incorporated in $\mathbf{A}$ as in~\cite{KleissPechsteinJuettlerTomar:2012}, since we only consider homogeneous Dirichlet conditions and by using the index sets $\mathbf{I}^{\left( j \right)}_0$ the corresponding entries in the local stiffness tensors $K^{\left( j \right)}$, $j = 1, \ldots, N_P$, are simply eliminated. We further note that $\boldsymbol{\lambda}$ is only unique up to an additive constant of $\ker \left( \mathbf{A}^\top \right)$. \\

As in~\cite{KleissPechsteinJuettlerTomar:2012}, the tensors $A^{\left( j \right)}_{\left(j, k\right)}$ and $A^{\left( k \right)}_{\left(j, k\right)}$ enforce $C^0$-continuity by linear constraints on the DoFs $y^{\left(j\right)}$ and $y^{\left(k\right)}$, which each are located on the interface $\Gamma^{\left(j, k\right)}$. For this, the tensor $A^{\left( j \right)}_{\left(j, k\right)}$ has to address the corresponding DoFs of patch $\Omega^{\left( j \right)}$ with indices in $\mathbf{I}_{\Gamma} (j, k)$. For that we make use of the tensor product structure of the solution space $V^{\left( j \right)}_h$. Each DoF $y^{\left( j \right)}_{\mathbf{i}}$, which represents the basis function $\hat{\beta}^{\left( j \right)}_\mathbf{i} \in V^{\left( j \right)}_h$, $\mathbf{i} = \left( i_1, i_2, i_3 \right) \in \mathbf{I}^{\left( j \right)}_0$, can be addressed by
\begin{equation}
\label{eq:address_DoF}
\begin{gathered}
y^{\left( j \right)}_{\mathbf{i}} = v_{\mathbf{i}} \cdot y^{\left( j \right)}, \\
v_{\mathbf{i}} = v_{i_1} \otimes v_{i_2} \otimes v_{i_3} \in \mathbb{R}^{\left(1, 1, 1\right) \times \left( \tilde{n}^{\left( j \right)}_1, \tilde{n}^{\left( j \right)}_2, \tilde{n}^{\left( j \right)}_3 \right)}, \\
v_{i_d} = \left[0, \ldots, 0, 1, 0, \ldots, 0 \right] \in \mathbb{R}^{1 \times \tilde{n}^{\left( j \right)}_d}, d = 1, 2, 3. \\
\hspace*{-2.40cm}\tikz[baseline=(arrow.base)]{
\node (arrow) {$i_d$-th entry};
\draw[->] (arrow.north) -- ++(0,0.5em);
}
\end{gathered}
\end{equation}
This justifies that for $\left(j, k\right) \in \mathcal{C}$ the tensors $A^{\left( m \right)}_{\left(j, k\right)}$, $m \in \{j, k\}$, have a rank-one representation, i.e.\ 
\begin{equation}
\label{eq:jump_tensor_rank-one}
A^{\left( m \right)}_{\left(j, k\right)} = A^{\left( m \right) \left( 1 \right)}_{\left(j, k\right)} \otimes A^{\left( m \right) \left( 2 \right)}_{\left(j, k\right)} \otimes A^{\left( m \right) \left( 3 \right)}_{\left(j, k\right)} \in \R^{\left( J^{\left(j, k\right)}_1, J^{\left(j, k\right)}_2, J^{\left(j, k\right)}_3 \right) \times \left( \tilde{n}^{\left( m \right)}_1, \tilde{n}^{\left( m \right)}_2, \tilde{n}^{\left( m \right)}_3\right)}
\end{equation}
where $A^{\left( m \right) \left( d \right)}_{\left(j, k\right)} \in \mathbb{R}^{J^{\left(j, k\right)}_d \times \tilde{n}^{\left( m \right)}_d}$, $J^{\left(j, k\right)}_d$ is the number of continuity constraints in dimension $d$, and the entries of these factor matrices depend on the underlying interface $\Gamma^{\left(j, k\right)}$. To simplify the notation, we will omit the index $\left(j, k\right)$ in the following. When transposing $\mathbf{A}$, the block structure of $\mathbf{A}$ is transposed as in the matrix case and in addition all factor matrices $A^{\left( m \right) \left( d \right)}$, $d \in \{1, 2, 3\}$, in~\eqref{eq:jump_tensor_rank-one} are transposed. \\

The $C^0$-continuity of the discrete approximation $y_h \colon \Omega \to \mathbb{R}$ can now be enforced by choosing suitable factor matrices in~\eqref{eq:jump_tensor_rank-one}. Their entries depend on which side of the unit cube the interface $\Gamma^{\left(j, k\right)}$ of the respective patch is located and the relationship between $y^{\left( j \right)}_{\mathbf{I}_{\Gamma} (j, k)}$ and $y^{\left( k \right)}_{\mathbf{I}_{\Gamma} (k, j)}$. By exploiting the tensor product structure of $V^{\left( m \right)}_h$, $m \in \{j, k\}$, we can split the $C^0$-continuity condition on the interface $\Gamma^{\left(j, k\right)}$ into $1$-dimensional $C^0$-continuity conditions and thus formulate the corresponding linear constraints in the factor matrix $A^{\left( m \right) \left( d \right)}$ for each dimension $d \in \{1, 2, 3\}$. \\

For $\left(j, k\right) \in \mathcal{C}$ let $\Omega^{\left( j \right)}$ and $\Omega^{\left( k \right)}$ be connected by the interface $\Gamma^{\left(j, k\right)}$ in dimension $d_1 \in \{ 1, 2, 3 \}$. Then $\mathbf{I}_{\Gamma} (j, k)$ and $\mathbf{I}_{\Gamma} (k, j)$ are $1$-dimensional in dimension $d_1$, which means that for each patch all basis functions on this interface $\hat{\beta}^{\left( j \right)}_{\mathbf{I}_{\Gamma} (j, k)}$, and thus the discrete function $y_h$, are constructed by only one univariate spline in dimension $d_1$, see~\eqref{eq:solution_tensor_product}. As this univariate spline must be located on one side of the unit cube for each patch in dimension $d_1$ and based on our assumption about the orientation of the patches, this is the first univariate spline for one patch and the last univariate spline for the other patch in dimension $d_1$. All basis functions on the interface $\Gamma^{\left(j, k\right)}$ of the corresponding patch have this one univariate spline in~\eqref{eq:solution_tensor_product}, and to ensure $C^0$-continuity on that interface, the DoFs of both patches that depend on this univariate spline must match in dimension $d_1$. From~\eqref{eq:address_DoF}, we can identify all DoFs of the basis functions in dimension $d_1$ for both patches $\Omega^{\left( j \right)}$ and $\Omega^{\left( k \right)}$ by using a row vector, i.e. $A^{\left( m \right) \left( d_1 \right)} \in \mathbb{R}^{1 \times \tilde{n}^{\left( m \right)}_{d_1}}$, $m \in \{j, k\}$. This row vector contains only zeros except for one $1$, which is either the first or the last entry of the vector for the corresponding patch, which depends on which of the two sides of the unit cube of the respective patch in dimension $d_1$ the interface $\Gamma^{\left(j, k\right)}$ is located. When we think of a dice again and for patch $\Omega^{\left( j \right)}$ the interface is located on side $1$ and for patch $\Omega^{\left( k \right)}$ on side $6$, then $d_1 = 1$ and we have
\begin{gather*}
A^{\left( j \right) \left( 1 \right)} = \left[ 1, 0, \ldots, 0 \right] \in \mathbb{R}^{1 \times \tilde{n}^{\left( j \right)}_{1}}, \\
A^{\left( k \right) \left( 1 \right)} = \left[ 0, \ldots, 0, 1 \right] \in \mathbb{R}^{1 \times \tilde{n}^{\left( k \right)}_{1}},
\end{gather*} 
as $d_1$-th factor matrix for $A^{\left( j \right)}$ and $A^{\left( k \right)}$ in~\eqref{eq:jump_tensor_rank-one}. This addresses the DoFs on the interface on the corresponding side of the respective patch in the dimension $d_1$ of the dice due to the tensor product structure. We note that we multiply one of the two factor matrices $A^{\left( j \right) \left( d_1 \right)}$ or $A^{\left( k \right) \left( d_1 \right)}$ by $-1$, since we want to subtract the resulting tensors from each other as in~\eqref{eq:DoFs_subtraction}.

The factor matrices in~\eqref{eq:jump_tensor_rank-one} for dimension $d_l \in \{1, 2, 3\} \setminus \{d_1\}$, $l \in \{2, 3\}$, depend on whether the patches $\Omega^{\left( j \right)}$ and $\Omega^{\left( k \right)}$ are \textit{fully matching} (cf.~\cite{KleissPechsteinJuettlerTomar:2012}) in this dimension $d_l$ on the interface $\Gamma^{\left(j, k\right)}$ or not. 

If the patches $\Omega^{\left( j \right)}$ and $\Omega^{\left( k \right)}$ are fully matching on $\Gamma^{\left(j, k\right)}$ in dimension $d_l \in \{1, 2, 3\} \setminus \{d_1\}$, $l \in \{2, 3\}$, then the knot vector $\xi^{\left( j \right)}_{d_l}$ is affinely related to the knot vector $\xi^{\left( k \right)}_{d_l}$ and the corresponding weights and degrees are equal. In the following, we will use the familiar term \textit{conforming}. With our assumption about the orientation of the patches, the two knot vectors actually match in this case, which in turn implies, that $\mathbb{S}^{\left( j \right)}_{\xi^{\left( j \right)}_{d_l}} = \mathbb{S}^{\left( k \right)}_{\xi^{\left( k \right)}_{d_l}}$ and $\tilde{n}^{\left( j \right)}_{d_l} = \tilde{n}^{\left( k \right)}_{d_l}$ holds. This means that the univariate factors of the basis functions of both patches $\Omega^{\left( j \right)}$ and $\Omega^{\left( k \right)}$ coincide in~\eqref{eq:solution_tensor_product} for this dimension $d_l$ on the interface $\Gamma^{\left(j, k\right)}$. To ensure $C^0$-continuity of the approximation, the $d_l$-th factor matrix in~\eqref{eq:jump_tensor_rank-one} for both patches must be a square Boolean matrix $A^{\left( m \right) \left( d_l \right)} \in \mathbb{R}^{\tilde{n}^{\left( m \right)}_{d_l}\times \tilde{n}^{\left( m \right)}_{d_l}}$, $m \in \{j, k\}$, whose rows are the vectors $v_{i_{d_{l}}} \in \mathbb{R}^{1 \times \tilde{n}^{\left( m \right)}_{d_l}}$ for dimension $d_l$ in~\eqref{eq:address_DoF}, which address the corresponding DoFs on the interface $\Gamma^{\left(j, k\right)}$. With our assumption about the orientation, these factor matrices in~\eqref{eq:jump_tensor_rank-one} are for both patches $\Omega^{\left( j \right)}$ and $\Omega^{\left( k \right)}$ the identity matrix, i.e. $A^{\left( j \right) \left( d_l \right)} = A^{\left( k \right) \left( d_l \right)} = \mathbb{I} \in \mathbb{R}^{\tilde{n}^{\left( j \right)}_{d_l} \times \tilde{n}^{\left( j \right)}_{d_l}}$. 

If the patches $\Omega^{\left( j \right)}$ and $\Omega^{\left( k \right)}$ are not fully matching or \textit{nonconforming} in dimension $d_l \in \{1, 2, 3\} \setminus \{d_1\}$, $l \in \{2, 3\}$, but the knot vector $\xi^{\left( j \right)}_{d_l}$ is obtained from the knot vector $\xi^{\left( k \right)}_{d_l}$ by one step of uniform h-refinement, then the univariate spline spaces $\mathbb{S}^{\left( j \right)}_{\xi^{\left( j \right)}_{d_l}}$ and $\mathbb{S}^{\left( k \right)}_{\xi^{\left( k \right)}_{d_l}}$ no longer coincide with each other and $\tilde{n}^{\left( j \right)}_{d_l} \neq \tilde{n}^{\left( k \right)}_{d_l}$. In~\cite{KleissPechsteinJuettlerTomar:2012} it was shown for the $2$-dimensional case, where the interface is a $1$-dimensional edge, that the DoFs of the finer patch $\Omega^{\left( j \right)}$ can be expressed as a linear combination of the DoFs of the coarser patch $\Omega^{\left( k \right)}$ on the interface $\Gamma^{\left(j, k\right)}$. We adopt this approach, which means that $\forall \left( i^{\left( j \right)}_{d_1}, i^{\left( j \right)}_{d_l}, i_{d_p} \right) \in \mathbf{I}_{\Gamma} (j, k)$
\begin{equation}
\label{eq:DOFs_linear_combination}
y^{\left( j \right)}_{\left( i^{\left( j \right)}_{d_1}, i^{\left( j \right)}_{d_l}, i_{d_p} \right)} = \sum^{\tilde{n}^{\left( k \right)}_{d_l}}_{i^{\left( k \right)}_{d_l} = 1} Z_{i^{\left( j \right)}_{d_l}, i^{\left(k\right)}_{d_l}} \, y^{\left( k \right)}_{\left( i^{\left( k \right)}_{d_1}, i^{\left( k \right)}_{d_l}, i_{d_p} \right)}, \quad \quad \left( i^{\left( k \right)}_{d_1}, i^{\left( k \right)}_{d_l}, i_{d_p} \right) \in \mathbf{I}_{\Gamma} (k, j),
\end{equation}
must hold to ensure $C^0$-continuity, where the linear coefficients $Z \in \mathbb{R}^{\tilde{n}^{\left( j \right)}_{d_l} \times \tilde{n}^{\left( k \right)}_{d_l}}$ can be obtained from the formula for h-refinement of B-spline basis functions (cf.~\cite[Section~5.3]{piegl}). We note that the index $i^{\left( m \right)}_{d_1} \in \{1, \tilde{n}^{\left( m \right)}_{d_1}\}$, $m \in \{j, k\}$, on both sides in~\eqref{eq:DOFs_linear_combination} is fixed and given by the location of the interface $\Gamma^{\left( j, k \right)}$ and for simplicity we assume that the patches $\Omega^{\left( j \right)}$ and $\Omega^{\left( k \right)}$ are conforming in dimension $d_p$, which is why we can assume that the index in dimension $d_p$ is the same. If the patches are nonconforming in dimension $d_p$ and the patch $\Omega^{\left( j \right)}$ is finer in this dimension, then this would result in a double sum in~\eqref{eq:DOFs_linear_combination}, so that the DoFs of the finer patch $\Omega^{\left( j \right)}$ can also be represented as a linear combination of the DoFs of the coarser patch $\Omega^{\left( k \right)}$. This means that the number of continuity constraints in dimension $d_l$ is given by the number of univariate splines of the finer patch $\Omega^{\left( j \right)}$, i.e. $A^{\left( j \right) \left( d_l \right)} \in \mathbb{R}^{\tilde{n}^{\left( j \right)}_{d_l} \times \tilde{n}^{\left( j \right)}_{d_l}}$ and $A^{\left( k \right) \left( d_l \right)} \in \mathbb{R}^{\tilde{n}^{\left( j \right)}_{d_l} \times \tilde{n}^{\left( k \right)}_{d_l}}$. The factor matrix $A^{\left( j \right) \left( d_l \right)}$ of the finer patch $\Omega^{\left( j \right)}$ is again a Boolean matrix, which addresses with its rows the corresponding DoFs on the interface $\Gamma^{\left(j, k \right)}$ as in~\eqref{eq:address_DoF}. With our assumption about the orientation of the patches, we can choose for that the identity matrix $\mathbb{I} \in \mathbb{R}^{\tilde{n}^{\left( j \right)}_{d_l} \times \tilde{n}^{\left( j \right)}_{d_l}}$. The factor matrix $A^{\left( k \right) \left( d_l \right)}$ of the coarser patch $\Omega^{\left( k \right)}$ is given by the coefficients matrix, i.e. $A^{\left( k \right) \left( d_l \right)} = Z$. 

\section{IETI-based low-rank method}
\label{section:IETI_method}

We now discuss how to compute an approximate solution of~\eqref{eq:SaddlePointFormulation} to obtain the linear coefficients $y^{\left( j \right)} \in \mathbb{R}^{\tilde{n}^{(j)}_1 \times \tilde{n}^{(j)}_2 \times \tilde{n}^{(j)}_3}$ for the patch-wise representation~\eqref{eq:discrete_function} of the discrete approximation $y_h \colon \Omega \to \mathbb{R}$ of~\eqref{eq:forward}. We consider this here using \textsc{MATLAB}. An extended block \amen\ method (implemented as \texttt{amen\_block\_solve.m} in the \textsc{TT-Toolbox}~\cite{tt-toolbox}), which allows us to solve large systems while preserving the block structure without assembling the whole equation system and returns the solution in a low-rank TT format (cf.~\cite{bdos-sb-2016},~\cite{ds-navier-2017}), could be used for solving~\eqref{eq:SaddlePointFormulation}, but our experiments have shown that this inevitably leads to the complication, that all blocks in $\mathbf{K}$ and $\mathbf{A}$ must have the same size for this solver. An alternative would be to fill the blocks that are too small with zeros in the corresponding entries, but this would greatly impair numerical stability and the performance of the method is not competitive. \\

As in~\cite{KleissPechsteinJuettlerTomar:2012}, our approach follows the idea of a FETI-like method by eliminating the primal variables $\mathbf{y}$ from the system~\eqref{eq:SaddlePointFormulation} and solving for the dual variables $\boldsymbol{\lambda}$. The primal variables $\mathbf{y}$ can then be easily recovered from the dual variables $\boldsymbol{\lambda}$. We achieve this by solving the Schur complement of~\eqref{eq:SaddlePointFormulation}, i.e.\ we search for $\boldsymbol{\lambda} = \left[ \lambda^{\left( 1 \right)}, \ldots, \lambda^{\left( \lvert \mathcal{C} \rvert \right)} \right]^\top$ that solves 
\begin{equation}
\label{eq:Schur_complement}
\mathbf{A} \, \mathbf{K}^{-1} \, \mathbf{A}^\top \; \boldsymbol{\lambda} = \mathbf{A} \, \mathbf{K}^{-1} \; \mathbf{f}.
\end{equation}
All blocks $\mathbf{A}$, $\mathbf{K}$, $\mathbf{f}$ and $\boldsymbol{\lambda}$ have the same structure as before and are in TT format. We note that tensors and tensor matrices in canonical format (such as~\eqref{eq:stiffness_tensor} and~\eqref{eq:jump_tensor_rank-one}) can easily be converted to TT format (cf.~\cite{Oseledets:2011}). Obviously, the number of variables in~\eqref{eq:Schur_complement} is smaller, since we only solve for the Lagrange multipliers $\boldsymbol{\lambda}$ and each block of $\boldsymbol{\lambda}$ represents an interface $\Gamma^{\left(j, k\right)}$ with $J^{\left(j, k\right)}_{d_2} J^{\left(j, k\right)}_{d_3}$ many linear constraints, where $J^{\left(j, k\right)}_{d_l} = \max \left( \left\{ \tilde{n}^{\left( j \right)}_{d_l}, \tilde{n}^{\left( k \right)}_{d_l} \right\} \right)$, $l \in \{ 2, 3 \}$, since $J^{\left(j, k\right)}_{d_1} = 1$.

We find an approximate solution $\boldsymbol{\lambda}^*$ of~\eqref{eq:Schur_complement} using a tensor block version of \ttgmres\ (implemented as \texttt{tt\_gmres\_block.m} in~\cite{tt-toolbox}, cf.~\cite{Dolgov:2013}). This method is in contrast to the block \amen\ method tensor-matrix free, i.e.\ we can define the linear operator $\mathbf{A} \, \mathbf{K}^{-1} \, \mathbf{A}^\top$ as a function handle and filling with zeros is not necessary. When solving~\eqref{eq:Schur_complement}, we take advantage of the fact that the blocks of $\boldsymbol{\lambda}$ can be represented as $2$-dimensional tensors in order to simplify the problem, since these tensors are only $1$-dimensional in the dimension $d_1$ in which the two corresponding patches are connected, i.e. $\lambda^{\left( \left(j, k\right) \right)} \in \mathbb{R}^{J^{\left(j, k\right)}_{d_1} \times J^{\left(j, k\right)}_{d_2} \times J^{\left(j, k\right)}_{d_3}} = \mathbb{R}^{1 \times J^{\left(j, k\right)}_{d_2} \times J^{\left(j, k\right)}_{d_3}} \cong \mathbb{R}^{J^{\left(j, k\right)}_{d_2} \times J^{\left(j, k\right)}_{d_3}}$. Therefore, the \ttgmres\ is applied to a $2$-dimensional linear block system in our setup, which reduces the complexity. When applying the function handle of the linear operator $\mathbf{A} \, \mathbf{K}^{-1} \, \mathbf{A}^\top$ on $\boldsymbol{\lambda}$, each block of $\boldsymbol{\lambda}$ is first reshaped into $3$-dimensional tensors using the \texttt{reshape.m} function of the \textsc{TT-Toolbox}, then the actual linear system is applied and finally all blocks of $\boldsymbol{\lambda}$ are reshaped back into $2$-dimensional tensors. When applying $\mathbf{K}^{-1}$ within the function handle, the patch-wise given linear systems are solved using the standard \amen\ method (implemented as \texttt{amen\_solve2.m} in~\cite{tt-toolbox}). \\

Although the system~\eqref{eq:Schur_complement} is smaller than the original system~\eqref{eq:SaddlePointFormulation}, preconditioners are still necessary to compute the variables $\boldsymbol{\lambda}$ in a reasonable number of iterations. This is because the system~\eqref{eq:Schur_complement} is generally ill-conditioned. We use a block diagonal tensor matrix with $\lvert \mathcal{C} \rvert$ many $3$-dimensional tensor matrices as diagonal blocks as a left preconditioner, i.e.\
\begin{equation}
\label{eq:preconditioned_Schu_complement}
\mathbf{P}^{-1} \, \mathbf{A} \, \mathbf{K}^{-1} \, \mathbf{A}^\top \; \boldsymbol{\lambda} = \mathbf{P}^{-1} \, \mathbf{A} \, \mathbf{K}^{-1} \; \mathbf{f},
\end{equation}
which we apply to the iterate $\boldsymbol{\lambda}^{\left( k \right)}$ after applying the linear operator $\mathbf{A} \, \mathbf{K}^{-1} \, \mathbf{A}^\top$ but before reshaping back to $2$-dimensional tensors in the function handle. 

For $\left(j, k\right) \in \mathcal{C}$ let $\Omega^{\left( j \right)}$ and $\Omega^{\left( k \right)}$ be connected by the interface $\Gamma^{\left(j, k\right)}$ in dimension $d_1 \in \{ 1, 2, 3 \}$ and assuming that the patch $\Omega^{\left( j \right)}$ is at least in one dimension $d_l \in \{1, 2, 3\} \setminus \{ d_1 \}$, $l \in \{2, 3\}$, finer than the patch $\Omega^{\left( k \right)}$. Then we explicitly set the diagonal blocks of $\mathbf{P}^{-1}$ with respect to the Lagrange multipliers of the interface $\Gamma^{\left(j, k\right)}$ as  
\begin{equation}
\label{eq:preconditioner_1}
P^{\left( \left(j, k\right) \right)} = A^{\left( j \right)}_{\left(j, k\right)} \cdot K^{\left( j \right)} \cdot {A^{\left( j \right)}_{\left(j, k\right)}}^{\mathrm{T}} \in \R^{\left( 1, J^{\left(j, k\right)}_{d_2}, J^{\left(j, k\right)}_{d_3} \right) \times \left( 1, J^{\left(j, k\right)}_{d_2}, J^{\left(j, k\right)}_{d_3} \right)},
\end{equation}
where $K^{\left( j \right)}$ is the low-rank stiffness tensor of patch $\Omega^{\left( j \right)}$ and $A^{\left( j \right)}_{\left(j, k\right)}$ is given by~\eqref{eq:jump_tensor_rank-one}, both in TT format. Here $\cdot$ denotes again a contracted product, namely the product of two TT matrices (cf.~\cite{Oseledets:2011}). We note that it is important that we use in~\eqref{eq:preconditioner_1} the information from the finer patch $\Omega^{\left( j \right)}$. In the conforming case, it has been shown that we can use the information from either patch and get similar results for the resulting two preconditioners.

\section{A PDE-constrained optimization model problem} 
\label{section:optimization}
We now want to discuss the discretization in both time and space of the optimization problem given on~\eqref{eq:optimization1} to~\eqref{eq:optimization4}, resulting in a large saddle point problem~\cite{saddlePoint, FEM}. Using an implicit Euler scheme for the time discretization of the PDE and the rectangle rule for the objective function leads to the time-discrete problem, which we then discretize in space using a Galerkin-based spatial discretization, which in turn leads to the discrete quadratic problem
\begin{align*}
\min_{y,u}&& \sum_{\ell=1}^{N_t}\frac{\tau}{2} \big ( (y_\ell-\hat{y}_\ell)^\top M (y_\ell-\hat{y}_t) &+ \alpha \, u_\ell^\top M u_\ell \big) & & \\
\mbox{s.t.}&& \frac{ M y_{\ell}- M y_{\ell-1}}{\tau} + Ky_{\ell} &= M u_{\ell} &&\mbox{ for } \ell = 1, \ldots, N_t,
\end{align*}
with the number of time steps $N_t$ corresponding to the time step size $\tau = T/{N_t}$. For the general case, $M$ and $K$ can be understood as mass or stiffness matrix of the corresponding geometry and all boundary conditions~\eqref{eq:optimization3} are incorporated in $M$ and $K$. The states are collected in a block vector $y = \left[ y_1, \ldots, y_{N_t} \right]^\top$ and similarly for the control $u$ and the desired state $\hat{y}$.

Such problems typically lead to saddle point systems as discussed in~\cite{saddlePoint, FEM}. We get to such a formulation by applying a Lagrangian formalism using a multiplier block vector $\mu = \left[ \mu_1, \ldots, \mu_{N_t} \right]^\top$ such that the \textit{Lagrangian} of the problem reads as
\begin{multline} 
\label{eq:objective_function}
\mathcal{L} \left( y, u, \mu \right) = \sum_{\ell=1}^{N_t} \biggl( \frac{\tau}{2} \Bigl( \bigl( y_\ell - \hat{y}_\ell \bigr)^\top M \bigl( y_\ell - \hat{y}_\ell \bigr) + \alpha \, u_\ell^\top M u_\ell \Bigr) \\
+ \mu_\ell^\top  \bigl( M y_{\ell} - M y_{\ell-1} + \tau K y_{\ell} - \tau M u_{\ell} \bigr) \biggr).
\end{multline}
Taking the derivative with respect to state $y$, control $u$ and Lagrange multiplier $\mu$ leads to the system 
\begin{equation} 
\label{eq:KKTSystem}
\begin{bmatrix} \tau \mathcal{M} & 0 & \mathcal{K}^\top \\ 0 &\tau \alpha \mathcal{M} & -\tau \mathcal{M} \\ \mathcal{K} & -\tau \mathcal{M} & 0 \end{bmatrix} \begin{bmatrix} y \\ u \\ \mu \end{bmatrix}  = \begin{bmatrix} \tau \mathcal{M}\hat{y} \\ 0 \\ 0 \end{bmatrix},
\end{equation}
where $\mathcal{M} = \mathbb{I} \otimes M$ and $\mathcal{K} = \mathbb{I} \otimes \tau K  + C \otimes M$, using the identity matrix $\mathbb{I} \in \mathbb{R}^{N_t\times N_t}$ and $C$ is representing the Euler scheme via
\begin{equation*}
C = \begin{bmatrix} 1 & 0 & 0 & \ldots & 0 \\ -1 & 1 & 0 & \ldots & 0 \\ 0 & -1 & 1 & \ldots & 0 \\ \vdots & & \ddots & \ddots&\\ 0 & \ldots & 0 & -1 & 1 \end{bmatrix}.
\end{equation*}
Note that in this derivation we used the same spline spaces for the state and control. It is also possible to have a different discretization for the control and this would make the system solver we use in the low-rank method more involved (cf.~\cite{BuengerDolgovStoll:2020}). The resulting equation system~\eqref{eq:KKTSystem} is a saddle point problem as described in~\cite{saddlePoint,stoll1,stoll2}.

We now define the linear system~\eqref{eq:KKTSystem} for multi-patch geometries, enforcing $C^0$-continuity for the state $y$ and the Lagrange multiplier $\mu,$ but not for the control $u$, since this is an algebraic variable (cf.~\cite{Herzog:2014}). We now include jump tensors as in the case of the elliptic problem. By introducing Larange multipliers with respect to the continuity constraints for $y$ and $\mu$ we obtain 
\begin{equation}
\label{eq:KKTSystem_2}
\begin{bmatrix} 
\tau \bar{M} & 0 & 0 & \bar{K}^\top & \bar{A} \\ 
0 & 0 & 0 & \bar{A}^\top & 0 \\ 
0 & 0 & \tau \alpha \bar{M} & -\tau \bar{M} & 0 \\ 
\bar{K} & \bar{A}^\top & -\tau \bar{M} & 0 & 0 \\ 
\bar{A} & 0 & 0 & 0 & 0 \\ 
\end{bmatrix} 
\begin{bmatrix} y \\ \lambda_y \\ u \\ \mu \\ \lambda_{\mu} \end{bmatrix}  = \begin{bmatrix} \tau \mathcal{M}\hat{y} \\ 0 \\ 0 \\ 0 \\ 0 \end{bmatrix},
\end{equation}
where $\bar{M}$ and $\bar{K}$ are block diagonal tensors as in~\eqref{eq:StiffnessMatrixMultiPatch} with diagonal blocks
\begin{equation*}
\bar{M}^{\left( j \right)} = \mathbb{I} \otimes M^{\left( j \right)}, \quad \bar{K}^{\left( j \right)} = \mathbb{I} \otimes \tau K^{\left( j \right)} + C \otimes M^{\left( j \right)}, \quad j = 1, \ldots, N_p,
\end{equation*}
where $M^{\left( j \right)}$ is defined by~\eqref{eq:massFinal} and $K^{\left( j \right)}$ by~\eqref{eq:stiffnessFinal}. Since the continuity should apply for all time steps, we set each block of this jump tensor $\bar{A}$ to $\mathbb{I} \otimes A^{\left( m \right)}_{\left(j, k\right)}$ where $A^{\left( m \right)}_{\left(j, k\right)}$, $\left( j, k \right) \in \mathcal{C}$, $m \in \left\{ j, k \right\}$, is defined as in section~\ref{section:jump_tensors}. The resulting saddle point problem~\eqref{eq:KKTSystem_2} typically becomes very large, depending on the number of time steps and refinement in the spatial discretization which is why we are solving the Schur complement for $u$
\begin{multline} 
\label{eq:OC_Schur_complement}
\left(\tau \alpha \bar{M} + \tau^3 
\begin{bmatrix}
\bar{M} & 0
\end{bmatrix}
\overbrace{{\begin{bmatrix}
\bar{K} & \bar{A}^\top \\
\bar{A} & 0
\end{bmatrix}}^{-\top}}^{= \bar{\mathcal{K}}^{-\top}}
\begin{bmatrix}
\bar{M} & 0 \\
0 & 0
\end{bmatrix}
\overbrace{{\begin{bmatrix}
\bar{K} & \bar{A}^\top \\
\bar{A} & 0
\end{bmatrix}}^{-1}}^{= \bar{\mathcal{K}}^{-1}}
\begin{bmatrix}
\bar{M} \\ 
0
\end{bmatrix}
\right) \, u = \\
\tau^2 
\begin{bmatrix}
\bar{M} & 0
\end{bmatrix}
{\begin{bmatrix}
\bar{K} & \bar{A}^\top \\
\bar{A} & 0
\end{bmatrix}}^{-\top}
\begin{bmatrix}
\bar{M} \\ 
0
\end{bmatrix} 
\hat{y},
\end{multline}
where the control is defined by $u = \left[ u^{\left( 1 \right)}, \ldots, u^{\left( N_P \right)} \right]^\top$. The blocks of $u$ are tensors of the form $u^{\left( j \right)} \in \mathbb{R}^{\tilde{n}^{(j)}_1 \times \tilde{n}^{(j)}_2 \times \tilde{n}^{(j)}_3 \times N_t}$ for each patch $\Omega^{\left( j \right)}$ and for all time steps. The state $y$ and multiplier $\mu$ can be easily recovered from the computed $u$. 

We consider solving~\eqref{eq:OC_Schur_complement} using \textsc{MATLAB}. In order to compute $u,$ we first transform al tensors in TT format and then use again the tensor block version of \ttgmres. The application of matrix vector product in~\eqref{eq:OC_Schur_complement} is computed by applying one operator after the other, which means that when applying~\eqref{eq:OC_Schur_complement} we have to solve a linear system once with $\bar{\mathcal{K}}$ and once with $\bar{\mathcal{K}}^{\top}$, for which we use the approach described in section~\ref{section:IETI_method} with the preconditioner described there. To reduce the number of iterations of \ttgmres, we apply the left preconditioner $\mathcal{P}$ whose block with respect to the patch $\Omega^{\left( j \right)}$ is defined as
\begin{equation}
\mathcal{P}^{\left( j \right)} = \tau \, \alpha \, \mathbb{I} \otimes \bigotimes^{3}_{d = 1} M^{\left( j \right) \left( d \right)}_{r_d}.
\end{equation}
Here $M^{\left( j \right) \left( d \right)}_{r_d}$ is one of the factor matrices in~\eqref{eq:stiffnessFinal} of $M^{\left( j \right)}$ and $r_d$ is chosen so that the norm of $M^{\left( j \right) \left( d \right)}_{r_d}$ is the largest for all $M^{\left( j \right) \left( d \right)}_{r}$, $r = 1, \ldots, R$, in this dimension $d$. 

Note that since we are using an iterative solver, to apply the operator of~\eqref{eq:OC_Schur_complement} in every iteration a flexible method such as {\sc FGMRES} will be ideally suited and we will tailor our approach in future research to this method as well designing a more sophisticated preconditioning strategy for this system.

\section{Numerical experiments} \label{section:numerics}

We now present the results of our numerical experiments. First, we investigate the error and performance of the method presented in section~\ref{section:IETI_method} with respect to different refinement levels and different solution tolerances. Secondly, we investigate the robustness of our in section~\ref{section:optimization} proposed method with respect to the penalty parameter $\alpha$. The experiments are conducted on two B-spline and one NURBS geometry with corresponding source functions or desired states for conforming and nonconforming patch discretizations. We point out that we restrict ourselves to low-rank multi-patch geometries, which means that the assembled mass and stiffness tensors of the individual patches have a low rank in the representation given by~\eqref{eq:massFinal} and~\eqref{eq:stiffnessFinal} (cf.~\cite{BuengerDolgovStoll:2020}). This choice of geometry allows us, at least for the B-spline geometries, to further reduce the cost of the assembly process. In more detail, we do not require a rich spline space for the interpolation of the weight functions as then refinement only takes place for the basis of the solution space. 

For our numerical experiments we used \textsc{MATLAB} R2022b on a desktop computer with an \textsc{AMD} Ryzen 5 5600X 6-core processor with 16 GB of RAM. Both geometries and specific functions from the \textsc{GeoPDEs 3.0} toolbox~\cite{geoPDEs} with the aid of the \textsc{NURBS Toolbox}~\cite{Spink:NURBS} were used. Computations in the TT-format were carried out using the \textsc{TT-Toolbox}~\cite{tt-toolbox}.

\subsection{Elliptic Problem}
\label{subsection:Elliptic_Problem}

Let $N_{DoFs}$ denote the total number of DoFs for the multi-patch geometry $\Omega$ and $N^{\left( j \right)}_{DoFs}$ the number of DoFs of patch $\Omega^{\left( j \right)}$ of $\Omega$. Thus, we define the relative $\mathcal{L}^2$-error to the analytical solution on a single patch and on the whole geometry as 
\begin{equation}
\label{eq:error_calculation}
\begin{gathered}
R^{\left( j \right)} \left( y_h, y_{sol} \right) = \frac{\lVert y_h|_{\Omega^{\left( j \right)}} - y_{sol}|_{\Omega^{\left( j \right)}} \rVert_{\mathcal{L}^2 \left( \Omega^{\left( j \right)} \right) }}{\lVert y_{sol}|_{\Omega^{\left( j \right)}} \rVert_{\mathcal{L}^2 \left( \Omega^{\left( j \right)} \right)}}, \\
R \left( y_h, y_{sol} \right) = \sum_{j = 1}^{N_p} \frac{N^{\left( j \right)}_{DoFs}}{N_{DoFs}} \, R^{\left( j \right)} \left( y_h, y_{sol} \right),
\end{gathered}
\end{equation}
where $y_h|_{\Omega^{\left( j \right)}}$ denotes the discrete approximation described by~\eqref{eq:discrete_function} for patch $\Omega^{\left( j \right)}$, $y_{sol} \colon \Omega \to \mathbb{R}$ is the analytical solution and $\lVert \cdot \rVert_{\mathcal{L}^2 \left( \Omega^{\left( j \right)} \right)}$ denotes the usual $\mathcal{L}^2$-norm. 

For each numerical experiment we measure the error to the analytical solution given by~\eqref{eq:error_calculation}, the number of iterations of the solver and the total time needed to compute the approximation, depending on different refinements of the solution space and different tolerances. Refinement in the conforming case is to be understood as starting with the original spline basis given by the geometry and then performing one step of uniform $h$-refinement with an increasing number of knots to be inserted between two existing knots in each dimension. In the nonconforming case, we do the same, but then undertake another one or two steps of uniform h-refinement for certain patches, in which only a single knot is inserted between two existing knots. Since we have to rely on functions of the toolbox \textsc{GeoPDEs 3.0}~\cite{geoPDEs} for the error calculation and for that we have to use the full coefficient vector, i.e. $\operatorname{vec} \left( y^{\left( j \right)} \right) \in \mathbb{R}^{\tilde{n}^{(j)}_1 \, \tilde{n}^{(j)}_2 \, \tilde{n}^{(j)}_3}$, this is currently our limiting factor for further increasing the number of DoFs. Here, $\varepsilon$ refers to the tolerance for solving the weight function interpolation system~\eqref{eq:weightfunction_linear_system} within \texttt{amen\_block\_solve.m} (cf.~\cite{BuengerDolgovStoll:2020}). For solving~\eqref{eq:preconditioned_Schu_complement} we use \texttt{tt\_gmres\_block.m}, where we set \texttt{max\_iters} to $10$, \texttt{restart} to $20$, and \texttt{tol} to $\varepsilon \cdot 10^2$. For solving the local linear systems defined by $K^{\left( j \right)}$ inside~\eqref{eq:preconditioned_Schu_complement} we use \texttt{amen\_solve2.m}\footnote{with parameters \texttt{nswp} $ = 20$, \texttt{kickrank} $= 2$ and \texttt{tol} $= \varepsilon \cdot 10$}. For the conforming cases, we compare the error and time of our method with the results of an approximation computed using \textsc{GeoPDEs 3.0}, which is limited to conforming geometries.  \\

\begin{figure}
\definecolor{mycolor1}{rgb}{0.00000,0.44700,0.74100}
\definecolor{mycolor2}{rgb}{0.85000,0.32500,0.09800}
\definecolor{mycolor3}{rgb}{0.92900,0.69400,0.12500}
\definecolor{mycolor4}{rgb}{0.49400,0.18400,0.55600}
\definecolor{mycolor5}{rgb}{0.46600,0.67400,0.18800}
\begin{subfigure}[t]{0.49\textwidth}
\centering
\includegraphics[width=\textwidth]{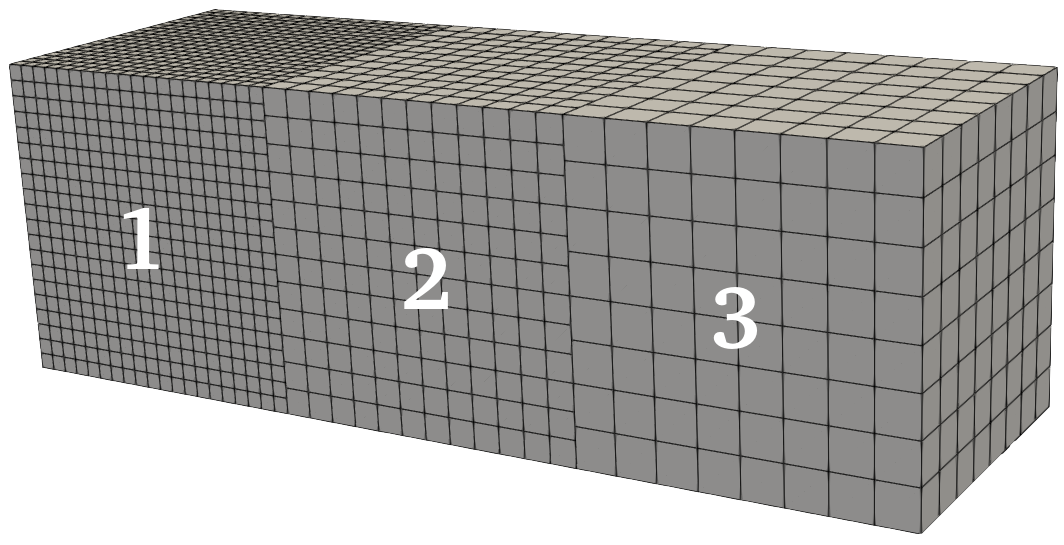}
\caption{3 cubes connected next to each other.}
\label{figure:3cubes_conf_1}
\end{subfigure}
\hfill
\begin{subfigure}[t]{0.49\textwidth}
\centering
\begin{tikzpicture}[font=\footnotesize]
\begin{axis}[
width=1.0\textwidth,
height=1.0\textwidth,
xmin=-50000,
xmax=700000,
xlabel={$N_{DoFs}$},
xtick = {0, 1e5, 2e5, 3e5, 4e5, 5e5, 6e5},
ymode=log,
ymin=1.0e-09,
ymax=10.0,
ylabel={$R ( y_{h}, y_{sol} )$},
yminorticks=true,
yminorgrids=true,
ymajorgrids=true,
xminorgrids=true,
ytick = {0.001, 1e-5, 1e-7},
extra y ticks = {1e-2, 1e-4, 1e-6, 1e-8},
extra y tick style={grid=none},
axis background/.style={fill=white}
]
\addplot [color=mycolor1, mark=*, mark size=1.0pt, mark options={solid, mycolor1}]
table[row sep=crcr]{
648	1.03105181947949\\
2187	0.0171357550857573\\
5184	0.0139116557160177\\
10125	0.00101485491762184\\
17496	0.000221418216395906\\
27783	0.00016848513229627\\
41472	0.000224686567165435\\
59049	0.000261418970422894\\
81000	0.000488683737259891\\
107811	0.00057232471057959\\
139968	0.000657352544190722\\
177957	0.000510659875102504\\
222264	0.000695312897705356\\
273375	0.000777785073663189\\
331776	0.000714929840806706\\
397953	0.000531433521186907\\
472392	0.000758136050689964\\
555579	0.000775406046240623\\
648000	0.000557151480081172\\
};
\addplot [color=mycolor2, mark=*, mark size=1.0pt, mark options={solid, mycolor2}]
table[row sep=crcr]{
648	1.03103388944605\\
2187	0.0171329355046082\\
5184	0.0139097985218257\\
10125	0.000986818862378438\\
17496	0.000122404326871156\\
27783	2.3120343250779e-05\\
41472	7.17644138781751e-06\\
59049	3.33600648051976e-06\\
81000	4.03821563524078e-06\\
107811	3.46431528333923e-06\\
139968	6.46305713476684e-06\\
177957	5.73183346629053e-06\\
222264	7.34182233990934e-06\\
273375	3.71354493526795e-06\\
331776	6.11091664612344e-06\\
397953	7.78922585672787e-06\\
472392	7.07944834390382e-06\\
555579	1.02146101731739e-05\\
648000	7.0269556188166e-06\\
};
\addplot [color=mycolor3, mark=*, mark size=1.0pt, mark options={solid, mycolor3}]
table[row sep=crcr]{
648	1.03103410525007\\
2187	0.0171328076872411\\
5184	0.0139097975528266\\
10125	0.000986806022476282\\
17496	0.000122387472433377\\
27783	2.29227809873791e-05\\
41472	6.41480253442148e-06\\
59049	2.2257129047289e-06\\
81000	9.03909641357906e-07\\
107811	4.1769344530865e-07\\
139968	2.09026416241021e-07\\
177957	1.28904463904627e-07\\
222264	9.18311137402375e-08\\
273375	1.00314615614909e-07\\
331776	9.2772115146761e-08\\
397953	6.3799943888134e-08\\
472392	7.54406799100933e-08\\
555579	6.8829430882887e-08\\
648000	5.77881469602793e-08\\
};
\addplot [color=mycolor5, mark=star, mark size=2.0pt, mark options={solid, mycolor5}]
table[row sep=crcr]{
648	0.800908323409017\\
2187	0.0171328211801889\\
5184	0.0139097975536947\\
10125	0.000986806028694988\\
17496	0.000122387456666969\\
27783	2.29227624533645e-05\\
41472	6.4147172690966e-06\\
59049	2.22525911857472e-06\\
};
\end{axis}
\end{tikzpicture}
\caption{Relative $\mathcal{L}^2$-error on $\Omega$.} \label{figure:3cubes_conf_2}
\end{subfigure}\\
\begin{subfigure}[t]{0.49\textwidth}
\centering
\begin{tikzpicture}[font=\footnotesize]

\begin{axis}[
width=1.0\textwidth,
height=1.0\textwidth,
xmin=-50000,
xmax=700000,
xlabel={$N_{DoFs}$},
xtick = {0, 1e5, 2e5, 3e5, 4e5, 5e5, 6e5},
ymin=0,
ymax=22,
ylabel={$\lvert \text{Iterations} \rvert$},
ytick = {5, 10, 15, 20},
axis background/.style={fill=white}
]
\addplot [color=mycolor1, mark=*, mark size=1.0pt, mark options={solid, mycolor1}]
table[row sep=crcr]{
648	1\\
2187	3\\
5184	4\\
10125	4\\
17496	5\\
27783	4\\
41472	5\\
59049	5\\
81000	5\\
107811	5\\
139968	5\\
177957	5\\
222264	4\\
273375	3\\
331776	4\\
397953	4\\
472392	5\\
555579	4\\
648000	4\\
};
\addplot [color=mycolor2, mark=*, mark size=1.0pt, mark options={solid, mycolor2}]
table[row sep=crcr]{
648	3\\
2187	6\\
5184	8\\
10125	9\\
17496	9\\
27783	8\\
41472	8\\
59049	9\\
81000	9\\
107811	11\\
139968	10\\
177957	10\\
222264	10\\
273375	12\\
331776	10\\
397953	10\\
472392	10\\
555579	11\\
648000	12\\
};
\addplot [color=mycolor3, mark=*, mark size=1.0pt, mark options={solid, mycolor3}]
table[row sep=crcr]{
648	3\\
2187	7\\
5184	11\\
10125	12\\
17496	13\\
27783	15\\
41472	15\\
59049	16\\
81000	14\\
107811	15\\
139968	16\\
177957	16\\
222264	16\\
273375	18\\
331776	14\\
397953	16\\
472392	17\\
555579	18\\
648000	15\\
};
\end{axis}
\end{tikzpicture}
\caption{Number of iterations.}\label{figure:3cubes_conf_3}
\end{subfigure}
\hfill
\begin{subfigure}[t]{0.49\textwidth}
\centering
\begin{tikzpicture}[font=\footnotesize]
\begin{axis}[
width=1.0\textwidth,
height=1.0\textwidth,
xmin=-50000,
xmax=700000,
xlabel={$N_{DoFs}$},
xtick = {0, 1e5, 2e5, 3e5, 4e5, 5e5, 6e5},
ymode=log,
ymin=0,
ymax=5000,
yminorticks=true,
ylabel={$s$},
axis background/.style={fill=white}
]
\addplot [color=mycolor1, mark=*, mark size=1.0pt, mark options={solid, mycolor1}]
table[row sep=crcr]{
648	0.818491\\
2187	3.584743\\
5184	1.301226\\
10125	1.636168\\
17496	5.504096\\
27783	5.742465\\
41472	6.852448\\
59049	9.198825\\
81000	13.284108\\
107811	20.651563\\
139968	23.060388\\
177957	43.647048\\
222264	31.141465\\
273375	35.023735\\
331776	45.125734\\
397953	50.665616\\
472392	62.034308\\
555579	71.337673\\
648000	89.40823\\
};
\addplot [color=mycolor2, mark=*, mark size=1.0pt, mark options={solid, mycolor2}]
table[row sep=crcr]{
648	0.968878\\
2187	4.248601\\
5184	2.219188\\
10125	3.497976\\
17496	6.896243\\
27783	6.616601\\
41472	7.976022\\
59049	10.6569\\
81000	14.610358\\
107811	23.496174\\
139968	26.209628\\
177957	45.958767\\
222264	33.883776\\
273375	42.792401\\
331776	48.48247\\
397953	54.430639\\
472392	67.076857\\
555579	76.498125\\
648000	95.103195\\
};
\addplot [color=mycolor3, mark=*, mark size=1.0pt, mark options={solid, mycolor3}]
table[row sep=crcr]{
648	0.98551\\
2187	5.809654\\
5184	4.090359\\
10125	8.021561\\
17496	16.673145\\
27783	21.651654\\
41472	23.938358\\
59049	29.976034\\
81000	29.67996\\
107811	41.92777\\
139968	45.964255\\
177957	72.955511\\
222264	61.9566\\
273375	77.519285\\
331776	75.845056\\
397953	88.111674\\
472392	105.027336\\
555579	111.428292\\
648000	127.765815\\
};
\addplot [color=mycolor5, mark=star, mark size=2.0pt, mark options={solid, mycolor5}]
table[row sep=crcr]{
648	0.28286\\
2187	6.086989\\
5184	32.043869\\
10125	93.772468\\
17496	211.571512\\
27783	388.388416\\
41472	654.978736\\
59049	1042.277362\\
};
\end{axis}
\end{tikzpicture}
\caption{Total time in seconds $s$.}\label{figure:3cubes_conf_4}
\end{subfigure}
\par\smallskip
\hspace{3em}
\begin{subfigure}[t]{0.9\textwidth}
\centering
\begin{tikzpicture}[font=\footnotesize]
\begin{axis}[
hide axis,
xmin=10,
xmax=50,
ymin=0,
ymax=0.4,
legend columns=5, 
legend style={draw=white!15!black,legend cell align=left, /tikz/every even column/.append style={column sep=0.4cm}}
]
\addlegendimage{mycolor1, mark=*, mark size=1.0pt, mark options={solid, mycolor1}}
\addlegendentry{$\varepsilon = 10^{-4}$};
\addlegendimage{mycolor2, mark=*, mark size=1.0pt, mark options={solid, mycolor2}}
\addlegendentry{$\varepsilon = 10^{-6}$};
\addlegendimage{mycolor3, mark=*, mark size=1.0pt, mark options={solid, mycolor3}}
\addlegendentry{$\varepsilon = 10^{-8}$};
\addlegendimage{color=mycolor5, mark=star, mark size=2.0pt, mark options={solid, mycolor5}}
\addlegendentry{\textsc{GeoPDEs 3.0}};
\end{axis}
\end{tikzpicture}
\end{subfigure}
\caption{Performance for different refinements and tolerances on the conforming multi-patch geometry with $3$ cubes shown in~\ref{figure:3cubes_conf_1} with~\eqref{eq:3cubes_sol} as analytical solution.} 
\label{figure:3cubes_conf_all} 
\end{figure}

\begin{figure}
\definecolor{mycolor1}{rgb}{0.00000,0.44700,0.74100}
\definecolor{mycolor2}{rgb}{0.85000,0.32500,0.09800}
\definecolor{mycolor3}{rgb}{0.92900,0.69400,0.12500}
\definecolor{mycolor4}{rgb}{0.49400,0.18400,0.55600}
\definecolor{mycolor5}{rgb}{0.46600,0.67400,0.18800}

\begin{subfigure}[t]{0.49\textwidth}
\centering
\begin{tikzpicture}[font=\footnotesize]
\begin{axis}[
width=1.0\textwidth,
height=1.0\textwidth,
xmin=-30000.0,
xmax=300000.0,
xlabel={$N_{DoFs}$},
xtick = {0, 5e4, 1e5, 1.5e5, 2e5, 2.5e5},
ymode=log,
ymin=1e-08,
ymax=10,
ylabel={$R ( y_{h}, y_{sol} )$},
yminorticks=true,
yminorgrids=true,
ymajorgrids=true,
xminorgrids=true,
ytick = {0.001, 1e-5, 1e-7},
extra y ticks = {1e-2, 1e-4, 1e-6},
extra y tick style={grid=none},
axis background/.style={fill=white}
]
\addplot [color=mycolor1, mark=*, mark size=1.0pt, mark options={solid, mycolor1}]
table[row sep=crcr]{
1288	0.205208258686187\\
3269	0.0118951779167648\\
6756	0.0107366615987753\\
12187	0.00118297437366724\\
20000	0.00243391708047078\\
30633	0.000778990842302327\\
44524	0.00117479293812995\\
62111	0.000418621948474882\\
83832	0.00041424164755153\\
110125	0.000728376919113332\\
141428	0.000778976593440084\\
178179	0.000615412644450334\\
220816	0.000795204611874767\\
269777	0.000728041798865737\\
};
\addplot [color=mycolor2, mark=*, mark size=1.0pt, mark options={solid, mycolor2}]
table[row sep=crcr]{
1288	0.205156572611812\\
3269	0.0118671148953167\\
6756	0.0107063649692326\\
12187	0.00113838902456158\\
20000	0.00215071343644465\\
30633	0.000258001525885474\\
44524	0.000551846394668969\\
62111	2.4674117117326e-05\\
83832	0.000130740502151756\\
110125	3.39898634813119e-05\\
141428	2.231858435593e-05\\
178179	1.11831719498114e-05\\
220816	6.12951517908267e-06\\
269777	1.08565232539129e-05\\
};
\addplot [color=mycolor3, mark=*, mark size=1.0pt, mark options={solid, mycolor3}]
table[row sep=crcr]{
1288	0.205156355298144\\
3269	0.0118671160138262\\
6756	0.0107063842176415\\
12187	0.00113833925922488\\
20000	0.00215004296660607\\
30633	0.000255104333397656\\
44524	0.000549952381530466\\
62111	2.24136127062433e-05\\
83832	0.000125752245040107\\
110125	3.09095766090551e-05\\
141428	1.71448810853472e-05\\
178179	6.65098420538688e-06\\
220816	3.36590398730743e-06\\
269777	1.69671341483533e-06\\
};
\end{axis}
\end{tikzpicture}
\caption{Relative $\mathcal{L}^2$-error on $\Omega$.} \label{figure:3cubes_non_conf_1}
\end{subfigure}
\hfill
\begin{subfigure}[t]{0.49\textwidth}
\centering
\begin{tikzpicture}[font=\footnotesize]
\begin{axis}[
width=1.0\textwidth,
height=1.0\textwidth,
xmin=-30000.0,
xmax=300000.0,
xlabel={$N_{DoFs}$},
xtick = {0, 5e4, 1e5, 1.5e5, 2e5, 2.5e5},
ymode=log,
ymin=1e-09,
ymax=10,
ylabel={$R^{ (j) } ( y_h, y_{sol} )$},
yminorticks=true,
yminorgrids=true,
ymajorgrids=true,
xminorgrids=true,
ytick = {1e-7},
extra y ticks = {1e-2, 1e-3, 1e-4, 1e-5, 1e-6, 1e-8},
extra y tick style={grid=none},
axis background/.style={fill=white},
]
\addplot [color=mycolor3, mark=triangle, mark size=2.0pt, mark options={solid, mycolor3}]
table[row sep=crcr]{
1288	0.017254892367609\\
3269	0.000532970737786049\\
6756	0.00023173851554713\\
12187	2.29240114887981e-05\\
20000	4.41715412645355e-06\\
30633	1.2026762022282e-06\\
44524	4.17705824025419e-07\\
62111	1.79427495234142e-07\\
83832	1.19081536450291e-07\\
110125	4.85969859305886e-08\\
141428	8.01675018545461e-08\\
178179	1.1639727314085e-07\\
220816	4.85606983954768e-08\\
269777	3.64687630186196e-08\\
};
\addplot [color=mycolor3, mark=+, mark size=2.0pt, mark options={solid, mycolor3}]
table[row sep=crcr]{
1288	0.128903094354364\\
3269	0.0171558568354474\\
6756	0.00502793086523259\\
12187	0.000532929555491034\\
20000	0.000989214545841403\\
30633	0.000231738754547774\\
44524	6.34160028934274e-05\\
62111	2.29240197320128e-05\\
83832	9.5349203757884e-06\\
110125	4.41721569186097e-06\\
141428	2.22620905028215e-06\\
178179	1.20119033245291e-06\\
220816	6.87745357041399e-07\\
269777	4.2975770476914e-07\\
};
\addplot [color=mycolor3, mark=o, mark size=2.0pt, mark options={solid, mycolor3}]
table[row sep=crcr]{
1288	0.960411146872572\\
3269	0.0732246236304393\\
6756	0.125979735284918\\
12187	0.017132777844118\\
20000	0.0395932422066809\\
30633	0.00499379900374233\\
44524	0.0139097791534829\\
62111	0.000532883027375153\\
83832	0.00379659023739663\\
110125	0.000986805845494499\\
141428	0.000578983433343774\\
178179	0.000231720131667015\\
220816	0.000122387471109688\\
269777	6.32763283246847e-05\\
};
\end{axis}

\end{tikzpicture}

\caption{Relative $\mathcal{L}^2$-error on $\Omega^{\left( j \right)}$.} \label{figure:3cubes_non_conf_2}
\end{subfigure} \\
\begin{subfigure}[t]{0.49\textwidth}
\centering
\begin{tikzpicture}[font=\footnotesize]

\begin{axis}[
width=1.0\textwidth,
height=1.0\textwidth,
xmin=-30000.0,
xmax=300000.0,
xlabel={$N_{DoFs}$},
xtick = {0, 5e4, 1e5, 1.5e5, 2e5, 2.5e5},
ymin=0,
ymax=40,
ylabel={$\lvert \text{Iterations} \rvert$},
ytick = {5, 10, 15, 20, 25, 30, 35},
axis background/.style={fill=white}
]
\addplot [color=mycolor1, mark=*, mark size=1.0pt, mark options={solid, mycolor1}]
table[row sep=crcr]{
1288	2\\
3269	4\\
6756	6\\
12187	7\\
20000	7\\
30633	5\\
44524	4\\
62111	5\\
83832	6\\
110125	5\\
141428	6\\
178179	6\\
220816	5\\
269777	4\\
};
\addplot [color=mycolor2, mark=*, mark size=1.0pt, mark options={solid, mycolor2}]
table[row sep=crcr]{
1288	4\\
3269	7\\
6756	8\\
12187	11\\
20000	12\\
30633	14\\
44524	14\\
62111	14\\
83832	16\\
110125	15\\
141428	13\\
178179	12\\
220816	12\\
269777	12\\
};
\addplot [color=mycolor3, mark=*, mark size=1.0pt, mark options={solid, mycolor3}]
table[row sep=crcr]{
1288	6\\
3269	9\\
6756	13\\
12187	15\\
20000	17\\
30633	19\\
44524	19\\
62111	24\\
83832	25\\
110125	26\\
141428	30\\
178179	31\\
220816	30\\
269777	32\\
};
\end{axis}
\end{tikzpicture}
\caption{Number of iterations.} \label{figure:3cubes_non_conf_3}
\end{subfigure}
\hfill
\begin{subfigure}[t]{0.49\textwidth}
\centering
\begin{tikzpicture}[font=\footnotesize]

\begin{axis}[
width=1.0\textwidth,
height=1.0\textwidth,
xmin=-30000.0,
xmax=300000.0,
xlabel={$N_{DoFs}$},
xtick = {0, 5e4, 1e5, 1.5e5, 2e5, 2.5e5},
ymode=log,
ymin=0,
ymax=500,
ylabel={$s$},
yminorticks=true,
axis background/.style={fill=white}
]
\addplot [color=mycolor1, mark=*, mark size=1.0pt, mark options={solid, mycolor1}]
table[row sep=crcr]{
1288	2.983547\\
3269	2.237678\\
6756	3.781795\\
12187	2.26528\\
20000	2.85822\\
30633	5.788589\\
44524	8.498343\\
62111	10.681473\\
83832	11.314467\\
110125	13.373102\\
141428	17.737334\\
178179	22.173873\\
220816	29.928799\\
269777	44.300543\\
};
\addplot [color=mycolor2, mark=*, mark size=1.0pt, mark options={solid, mycolor2}]
table[row sep=crcr]{
1288	3.015595\\
3269	2.663801\\
6756	4.293339\\
12187	3.066954\\
20000	3.668007\\
30633	7.655791\\
44524	10.293674\\
62111	13.91474\\
83832	16.40048\\
110125	17.262569\\
141428	20.224797\\
178179	25.260717\\
220816	33.062583\\
269777	44.774032\\
};
\addplot [color=mycolor3, mark=*, mark size=1.0pt, mark options={solid, mycolor3}]
table[row sep=crcr]{
1288	3.412105\\
3269	3.811806\\
6756	7.228192\\
12187	7.442318\\
20000	9.988555\\
30633	15.233851\\
44524	17.095504\\
62111	28.089748\\
83832	33.094574\\
110125	40.768398\\
141428	51.740414\\
178179	65.719172\\
220816	71.65612\\
269777	98.322773\\
};
\end{axis}

\end{tikzpicture}
\caption{Total time in seconds $s$.} \label{figure:3cubes_non_conf_4}
\end{subfigure}
\par\smallskip
\hspace{3em}
\centering
\begin{subfigure}[t]{0.9\textwidth}
\centering
\begin{tikzpicture}[font=\footnotesize]
\begin{axis}[
hide axis,
xmin=10,
xmax=50,
ymin=0,
ymax=0.4,
legend columns=3, 
legend style={draw=white!15!black,legend cell align=left, /tikz/every even column/.append style={column sep=0.4cm}}
]
\addlegendimage{mycolor1, mark=*, mark size=1.0pt, mark options={solid, mycolor1}}
\addlegendentry{$\varepsilon = 10^{-4}$};
\addlegendimage{mycolor2, mark=*, mark size=1.0pt, mark options={solid, mycolor2}}
\addlegendentry{$\varepsilon = 10^{-6}$};
\addlegendimage{mycolor3, forget plot}
\addlegendentry{$\varepsilon = 10^{-8}$};
\addlegendimage{mycolor3, mark=triangle, mark size=2.0pt, mark options={solid, mycolor3}}
\addlegendentry{patch 1};
\addlegendimage{mycolor3, mark=+, mark size=2.0pt, mark options={solid, mycolor3}}
\addlegendentry{patch 2};
\addlegendimage{mycolor3, mark=o, mark size=2.0pt, mark options={solid, mycolor3}}
\addlegendentry{patch 3};
\end{axis}
\end{tikzpicture}
\end{subfigure}
\caption{Performance for different refinements and tolerances on the nonconforming multi-patch geometry with $3$ cubes shown in~\ref{figure:3cubes_conf_1} with~\eqref{eq:3cubes_sol} as analytical solution.} \label{figure:3cubes_non_conf_all} 
\end{figure}

We first study the multi-patch B-spline geometry shown in Figure \ref{figure:3cubes_conf_1} where the analytical solution of~\eqref{eq:forward} for this experiment is given by
\begin{equation}
\label{eq:3cubes_sol}
y_{sol} \left( x, y, z \right) = \sin \left( \sin\left( y \, \pi \right) \, \sin\left(z \, 2 \pi \right) \, \sin\left(x \, \pi \right)\right) \, \sin\left(y \, \pi \right) \, \sin\left(z \, 2 \pi \right) \, \sin\left(x \, \pi \right).
\end{equation}
We use B-splines of degree $5$ in both the conforming case and the nonconforming case. In the conforming case, we start with $N^{\left( j \right)}_{DoFs} = 216$ and we increase the number of knots to be inserted until we reach $N^{\left( j \right)}_{DoFs} = 216.000$ DoFs for $j = 1, 2, 3$. In the nonconforming case, we start with the division $N^{\left( 1 \right)}_{DoFs} = 729$, $N^{\left( 2 \right)}_{DoFs} = 343$, $N^{\left( 3 \right)}_{DoFs} = 216$ and end with $N^{\left( 1 \right)}_{DoFs} = 226.981$, $N^{\left( 2 \right)}_{DoFs} = 35.937$, $N^{\left( 3 \right)}_{DoFs} = 6.859$. 

Figure \ref{figure:3cubes_conf_2} shows the relative $\mathcal{L}^2$-error $R \left( y_h, y_{sol} \right)$ on the whole multi-patch geometry $\Omega$ for different $\varepsilon$ and depending on $N_{DoFs}$. We see at a certain point, further refinement no longer reduces the error for a given fixed tolerance $\varepsilon$, since for $\varepsilon = 10^{-4}$ the error $R \left( y_h, y_{sol} \right)$ does not decrease further after the sixth step of increasing $N_{DoFs}$, but increases and approaches $\varepsilon = 10^{-3}$. Similar can be recognised for $\varepsilon = 10^{-6}, 10^{-8}$. We can derive from that, that if the error should be reduced by refinement, this must be done together with a calibration of the tolerance $\varepsilon$. We can also see in Figure \ref{figure:3cubes_conf_2} the error of the approximation generated using \textsc{GeoPDEs 3.0}, where we are only able to compute the solutions until the eighth step of increasing $N_{DoFs}$, due to higher memory requirements. Nevertheless, our method and \textsc{GeoPDEs 3.0} show high agreement. We can derive from the number of iterations shown in Figure \ref{figure:3cubes_conf_3} that they depend mainly on $\varepsilon$ and we see in Figure \ref{figure:3cubes_conf_4} that the timing for our method shows a very benign growth with increasing number of $N_{DoFs}$. In summary, it can be deduced from the results that the error behaviour of the method presented in section~\ref{section:IETI_method} shows very mild dependence on the number of degrees of freedom but that also a higher accuracy requires further adjustment of the tolerance levels.

Figure \ref{figure:3cubes_non_conf_all} shows a similar set of results for the nonconforming case. The decrease of the error $R \left( y_h, y_{sol} \right)$ is slower than in the conforming case, which is due to the fact that the patches are refined differently and therefore error decrease at different rates. This is well illustrated in Figure \ref{figure:3cubes_non_conf_2}. Here the error $R^{\left( j \right)} \left( y_h, y_{sol} \right)$ is shown for the tolerance $\varepsilon = 10^{-8}$ on a single patch $\Omega^{\left( j \right)}$ depending on $N_{DoFs}$. In our numerical tests it has been shown that an actual difference in the error $R^{\left( j \right)} \left( y_h, y_{sol} \right)$ for the patches with different refinement is only recognizable for higher tolerances. Here too, the error behaviour is as expected. The patch with the highest refinement, patch 1 in Figure \ref{figure:3cubes_conf_1}, converges at the fastest rate to an error of $10^{-7}$. It is to be expected that the error for the other two patches will correspond to this if further refinements are carried out. The number of iterations shown in Figure \ref{figure:3cubes_non_conf_3} is higher than in the conforming case, but this is due to the fact that the system~\eqref{eq:preconditioned_Schu_complement} is more complex, since the factor matrices in~\eqref{eq:jump_tensor_rank-one} are not just Boolean matrices in that case, but have more structure. Nevertheless, the method converges in a satisfying time also in the nonconforming case. \\

\begin{figure}
\definecolor{mycolor1}{rgb}{0.00000,0.44700,0.74100}
\definecolor{mycolor2}{rgb}{0.85000,0.32500,0.09800}
\definecolor{mycolor3}{rgb}{0.92900,0.69400,0.12500}
\definecolor{mycolor4}{rgb}{0.49400,0.18400,0.55600}
\definecolor{mycolor5}{rgb}{0.46600,0.67400,0.18800}

\begin{subfigure}[t]{0.49\textwidth}
\centering
\includegraphics[width=0.75\textwidth]{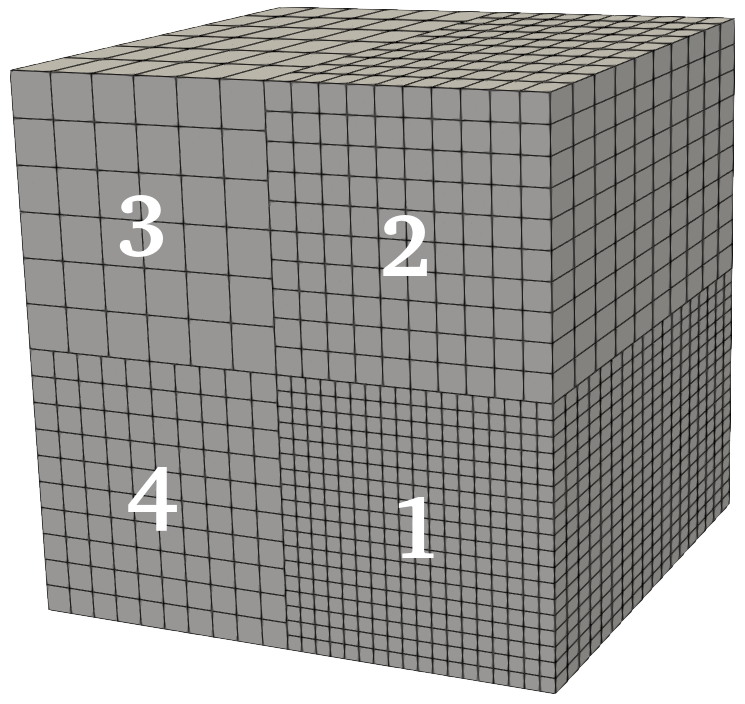}
\caption{4 cuboids connected to form a cube.}
\label{figure:4cubes_conf_1}
\end{subfigure}
\hfill
\begin{subfigure}[t]{0.49\textwidth}
\centering
\begin{tikzpicture}[font=\footnotesize]
\begin{axis}[
width=1.0\textwidth,
height=1.0\textwidth,
xmin=-100000,
xmax=1650000,
xlabel={$N_{DoFs}$},
ymode=log,
ymin=1.0e-08,
ymax=10,
yminorticks=true,
ylabel={$R ( y_{h}, y_{sol} )$},
yminorticks=true,
yminorgrids=true,
ymajorgrids=true,
xminorgrids=true,
ytick = {1e-4, 1e-6},
extra y ticks = {1e-2, 1e-3, 1e-5, 1e-7},
extra y tick style={grid=none},
axis background/.style={fill=white}
]
\addplot [color=mycolor1, mark=*, mark size=1.0pt, mark options={solid, mycolor1}]
table[row sep=crcr]{
256	1.00328625925385\\
2048	0.0311416714868417\\
6912	0.00165041124112434\\
16384	0.000312883220952315\\
32000	9.9838079551052e-05\\
55296	5.14232726443096e-05\\
87808	6.3915316413953e-05\\
131072	4.27139204682297e-05\\
186624	6.57619604020081e-05\\
256000	3.14420629172401e-05\\
340736	6.13261734776416e-05\\
442368	5.16987011797469e-05\\
562432	5.67246420654621e-05\\
702464	5.79875880473829e-05\\
864000	7.03903949914543e-05\\
1048576	7.62618048577874e-05\\
1257728	9.47565390151114e-05\\
1492992	7.20231315628152e-05\\
};
\addplot [color=mycolor2, mark=*, mark size=1.0pt, mark options={solid, mycolor2}]
table[row sep=crcr]{
256	1.00328625925384\\
2048	0.031141673966846\\
6912	0.00165039984537039\\
16384	0.000311793834861233\\
32000	9.79732369634889e-05\\
55296	4.02819702473102e-05\\
87808	1.95622429211709e-05\\
131072	1.06425573125217e-05\\
186624	6.29628786865847e-06\\
256000	3.98562076956595e-06\\
340736	2.65786381066459e-06\\
442368	1.87984163495194e-06\\
562432	1.41714018240327e-06\\
702464	1.12467982552409e-06\\
864000	9.401236929286e-07\\
1048576	8.69293936792985e-07\\
1257728	8.46109820497063e-07\\
1492992	1.23715761879672e-06\\
};
\addplot [color=mycolor3, mark=*, mark size=1.0pt, mark options={solid, mycolor3}]
table[row sep=crcr]{
256	1.00328625925342\\
2048	0.0311416740435769\\
6912	0.00165039986009687\\
16384	0.000311795172298005\\
32000	9.79731767747129e-05\\
55296	4.02819624624305e-05\\
87808	1.95567089802648e-05\\
131072	1.06317416048095e-05\\
186624	6.27429917462768e-06\\
256000	3.94106792875189e-06\\
340736	2.59992744809583e-06\\
442368	1.78446681926225e-06\\
562432	1.26545339931621e-06\\
702464	9.22332774753145e-07\\
864000	6.8810987666212e-07\\
1048576	5.23827890884015e-07\\
1257728	4.0588018130965e-07\\
1492992	3.19281220890732e-07\\
};
\addplot [color=mycolor5, mark=star, mark size=2.0pt, mark options={solid, mycolor5}]
table[row sep=crcr]{
256	1.0015501204792\\
2048	0.0311416886654154\\
6912	0.00165039991110988\\
16384	0.000311795175745568\\
32000	9.79731771394691e-05\\
55296	4.02819622751586e-05\\
87808	1.95567085429363e-05\\
};
\end{axis}
\end{tikzpicture}
\caption{Relative $\mathcal{L}^2$-error on $\Omega$.} \label{figure:4cubes_conf_2}
\end{subfigure} \\
\begin{subfigure}[t]{0.49\textwidth}
\centering
\begin{tikzpicture}[font=\footnotesize]

\begin{axis}[
width=1.0\textwidth,
height=1.0\textwidth,
xmin=-100000,
xmax=1650000,
xlabel={$N_{DoFs}$},
ymin=0,
ymax=60,
ytick={5, 10, 15, 20, 25, 30, 35, 40, 45, 50, 55},
ylabel={$\lvert \text{Iterations} \rvert$},
axis background/.style={fill=white}
]
\addplot [color=mycolor1, mark=*, mark size=1.0pt, mark options={solid, mycolor1}]
table[row sep=crcr]{
256	5\\
2048	5\\
6912	5\\
16384	9\\
32000	9\\
55296	9\\
87808	25\\
131072	26\\
186624	25\\
256000	8\\
340736	10\\
442368	10\\
562432	12\\
702464	8\\
864000	9\\
1048576	25\\
1257728	7\\
1492992	9\\
};
\addplot [color=mycolor2, mark=*, mark size=1.0pt, mark options={solid, mycolor2}]
table[row sep=crcr]{
256	9\\
2048	12\\
6912	13\\
16384	18\\
32000	17\\
55296	19\\
87808	19\\
131072	25\\
186624	25\\
256000	25\\
340736	29\\
442368	29\\
562432	33\\
702464	31\\
864000	29\\
1048576	28\\
1257728	27\\
1492992	34\\
};
\addplot [color=mycolor3, mark=*, mark size=1.0pt, mark options={solid, mycolor3}]
table[row sep=crcr]{
256	12\\
2048	18\\
6912	21\\
16384	27\\
32000	29\\
55296	33\\
87808	35\\
131072	36\\
186624	37\\
256000	41\\
340736	40\\
442368	40\\
562432	41\\
702464	42\\
864000	48\\
1048576	47\\
1257728	49\\
1492992	50\\
};
\end{axis}
\end{tikzpicture}
\caption{Number of iterations.}\label{figure:4cubes_conf_3}
\end{subfigure}
\hfill
\begin{subfigure}[t]{0.49\textwidth}
\centering
\begin{tikzpicture}[font=\footnotesize]

\begin{axis}[
width=1.0\textwidth,
height=1.0\textwidth,
xmin=-100000,
xmax=1650000,
xlabel={$N_{DoFs}$},
ymode=log,
ymin=0,
ymax=2000,
yminorticks=true,
ylabel={$s$},
axis background/.style={fill=white}
]
\addplot [color=mycolor1, mark=*, mark size=1.0pt, mark options={solid, mycolor1}]
table[row sep=crcr]{
256	0.820522\\
2048	0.519859\\
6912	0.670286\\
16384	1.224692\\
32000	1.672305\\
55296	2.354442\\
87808	8.513765\\
131072	10.396834\\
186624	12.374312\\
256000	8.690275\\
340736	11.944473\\
442368	15.234664\\
562432	19.311472\\
702464	22.361728\\
864000	27.986214\\
1048576	45.577286\\
1257728	39.82229\\
1492992	49.41267\\
};
\addplot [color=mycolor2, mark=*, mark size=1.0pt, mark options={solid, mycolor2}]
table[row sep=crcr]{
256	0.500931\\
2048	0.913637\\
6912	1.477724\\
16384	2.814316\\
32000	3.495304\\
55296	5.376186\\
87808	13.196025\\
131072	23.182926\\
186624	29.419596\\
256000	36.334446\\
340736	52.430504\\
442368	57.112057\\
562432	75.571453\\
702464	79.52087\\
864000	82.341638\\
1048576	93.273587\\
1257728	101.629043\\
1492992	139.645744\\
};
\addplot [color=mycolor3, mark=*, mark size=1.0pt, mark options={solid, mycolor3}]
table[row sep=crcr]{
256	0.547328\\
2048	1.453547\\
6912	2.934872\\
16384	6.37858\\
32000	11.206195\\
55296	18.776982\\
87808	38.256158\\
131072	53.909676\\
186624	70.850939\\
256000	102.378932\\
340736	122.717203\\
442368	141.840411\\
562432	165.460624\\
702464	187.205412\\
864000	241.188608\\
1048576	252.988009\\
1257728	276.18508\\
1492992	309.625477\\
};
\addplot [color=mycolor5, mark=star, mark size=2.0pt, mark options={solid, mycolor5}]
table[row sep=crcr]{
256	2.246144\\
2048	1.921636\\
6912	8.747327\\
16384	24.550477\\
32000	54.494688\\
55296	109.02141\\
87808	191.748484\\
};
\end{axis}            
\end{tikzpicture}
\caption{Total time in seconds $s$.} \label{figure:4cubes_conf_4}
\end{subfigure}
\par\smallskip
\hspace{3em}
\begin{subfigure}[t]{0.9\textwidth}
\centering
\begin{tikzpicture}[font=\footnotesize]
\begin{axis}[
hide axis,
xmin=10,
xmax=50,
ymin=0,
ymax=0.4,
legend columns=5, 
legend style={draw=white!15!black,legend cell align=left, /tikz/every even column/.append style={column sep=0.4cm}}
]
\addlegendimage{mycolor1, mark=*, mark size=1.0pt, mark options={solid, mycolor1}}
\addlegendentry{$\varepsilon = 10^{-4}$};
\addlegendimage{mycolor2, mark=*, mark size=1.0pt, mark options={solid, mycolor2}}
\addlegendentry{$\varepsilon = 10^{-6}$};
\addlegendimage{mycolor3, mark=*, mark size=1.0pt, mark options={solid, mycolor3}}
\addlegendentry{$\varepsilon = 10^{-8}$};
\addlegendimage{color=mycolor5, mark=star, mark size=2.0pt, mark options={solid, mycolor5}}
\addlegendentry{\textsc{GeoPDEs 3.0}};
\end{axis}
\end{tikzpicture}
\end{subfigure}
\caption{Performance for different refinements and tolerances on the conforming multi-patch geometry with $4$ cuboids shown in~\ref{figure:4cubes_conf_1} with~\eqref{eq:4cubes_sol} as analytical solution.} 
\label{figure:4cubes_conf_all} 
\end{figure}

The second B-spline geometry consisting of four cuboids connected to a cube is shown in Figure \ref{figure:4cubes_conf_1}, where the analytical solution is given by
\begin{equation}
\label{eq:4cubes_sol}
y_{sol} \left( x, y, z \right) = \sin\left(x \, 3 \pi \right) \, \sin \left(y \, \pi \right) \, \sin \left( z \, \pi \right).
\end{equation}
We use B-splines of degree $3$ in both the conforming case (results in Figure \ref{figure:4cubes_conf_all}) and the nonconforming case (results in Figure \ref{figure:4cubes_non_conf_all}). In the conforming case, we start with $N^{\left( j \right)}_{DoFs} = 64$ and end with $N^{\left( j \right)}_{DoFs} = 373.248$ for $j = 1, \ldots, 4$. In the nonconforming case, we start with $N^{\left( 1 \right)}_{DoFs} = 343$, $N^{\left( 2 \right)}_{DoFs}, N^{\left( 4 \right)}_{DoFs} = 125$, $N^{\left( 3 \right)}_{DoFs} = 64$ and end with $N^{\left( 1 \right)}_{DoFs} = 357.911$, $N^{\left( 2 \right)}_{DoFs}, N^{\left( 4 \right)}_{DoFs} = 50.653$, $N^{\left( 3 \right)}_{DoFs} = 8000$. We observe similar trends for this domain. It is clear that the iteration numbers do increase for the nonconforming domain and that the tolerance for an increased number of DoFs also needs further adjustment of the tolerances as stagnation of accuracy can be observed when the tolerances are not decreased accordingly. \\

\begin{figure}
\definecolor{mycolor1}{rgb}{0.00000,0.44700,0.74100}
\definecolor{mycolor2}{rgb}{0.85000,0.32500,0.09800}
\definecolor{mycolor3}{rgb}{0.92900,0.69400,0.12500}
\definecolor{mycolor4}{rgb}{0.49400,0.18400,0.55600}
\definecolor{mycolor5}{rgb}{0.46600,0.67400,0.18800}

\begin{subfigure}[t]{0.49\textwidth}
\centering
\begin{tikzpicture}[font=\footnotesize]
\begin{axis}[
width=1.0\textwidth,
height=1.0\textwidth,
xmin=-30000.0,
xmax=500000.0,
xlabel={$N_{DoFs}$},
ymode=log,
ymin=1.0e-07,
ymax=10,
ylabel={$R ( y_{h}, y_{sol} )$},
yminorticks=true,
yminorgrids=true,
ymajorgrids=true,
xminorgrids=true,
ytick = {1e-4, 1e-6},
extra y ticks = {1e-2, 1e-3, 1e-5},
extra y tick style={grid=none},
axis background/.style={fill=white}
]
\addplot [color=mycolor1, mark=*, mark size=1.0pt, mark options={solid, mycolor1}]
table[row sep=crcr]{
657	0.388285279585054\\
2142	0.042871317188034\\
5049	0.0666892253110478\\
9864	0.00666310804852494\\
17073	0.00524303118041161\\
27162	0.00165336024912346\\
40617	0.000658436882109003\\
57924	0.000314694018090945\\
79569	0.000167257417103516\\
106038	0.000115332158580764\\
137817	6.56389048299174e-05\\
175392	6.61678500053947e-05\\
219249	4.42194313275574e-05\\
269874	4.85128758328267e-05\\
327753	6.53523793665467e-05\\
393372	2.99079626887028e-05\\
467217	8.84655415108044e-05\\
};
\addplot [color=mycolor2, mark=*, mark size=1.0pt, mark options={solid, mycolor2}]
table[row sep=crcr]{
657	0.388523721015555\\
2142	0.043084376891088\\
5049	0.0636274463596355\\
9864	0.00670673354113813\\
17073	0.00535643571282472\\
27162	0.00167946212666223\\
40617	0.000678128537460175\\
57924	0.000313546607009634\\
79569	0.000163116243407696\\
106038	9.22631296325271e-05\\
137817	5.68121753065402e-05\\
175392	3.56487714567425e-05\\
219249	2.37907261054871e-05\\
269874	1.64733207498933e-05\\
327753	1.17819259138602e-05\\
393372	8.85158468142813e-06\\
467217	7.45969506425613e-06\\
};
\addplot [color=mycolor3, mark=*, mark size=1.0pt, mark options={solid, mycolor3}]
table[row sep=crcr]{
657	0.388523936391866\\
2142	0.0430844046897564\\
5049	0.0636270060221259\\
9864	0.00670691897294138\\
17073	0.0053566295639672\\
27162	0.00167962747901974\\
40617	0.000678239870103585\\
57924	0.000313520972587005\\
79569	0.000163156588289496\\
106038	9.22704318497091e-05\\
137817	5.58234721385346e-05\\
175392	3.56038509290415e-05\\
219249	2.3715052626376e-05\\
269874	1.63715744594156e-05\\
327753	1.16480127720858e-05\\
393372	8.50258385885594e-06\\
467217	6.34516824648212e-06\\
};
\end{axis}
\end{tikzpicture}
\caption{Relative $\mathcal{L}^2$-error on $\Omega$.} \label{figure:4cubes_non_conf_1}
\end{subfigure}
\hfill
\begin{subfigure}[t]{0.49\textwidth}
\centering
\begin{tikzpicture}[font=\footnotesize]
\begin{axis}[
width=1.0\textwidth,
height=1.0\textwidth,
xmin=-30000.0,
xmax=500000.0,
xlabel={$N_{DoFs}$},
ymode=log,
ymin=1e-07,
ymax=10,
ylabel={$R^{ (j) } ( y_h, y_{sol} )$},
yminorticks=true,
yminorgrids=true,
ymajorgrids=true,
xminorgrids=true,
ytick = {1e-6},
extra y ticks = {1e-2, 1e-3, 1e-4, 1e-5},
extra y tick style={grid=none},
axis background/.style={fill=white},
]
\addplot [color=mycolor3, mark=triangle, mark size=2.0pt, mark options={solid, mycolor3}]
table[row sep=crcr]{
657	0.214503454518254\\
2142	0.0226690602187714\\
5049	0.0343940060841312\\
9864	0.00324160680399295\\
17073	0.00273673000380782\\
27162	0.000777538719075499\\
40617	0.000293886774666138\\
57924	0.000127810298623085\\
79569	6.32788516970012e-05\\
106038	3.42147765929802e-05\\
137817	1.98998045555826e-05\\
175392	1.22467457387221e-05\\
219249	7.89775843514097e-06\\
269874	5.2928940983366e-06\\
327753	3.66434989380184e-06\\
393372	2.60791742384215e-06\\
467217	1.9007531246354e-06\\
};
\addplot [color=mycolor3, mark=+, mark size=2.0pt, mark options={solid, mycolor3}]
table[row sep=crcr]{
657	0.470315924070121\\
2142	0.0624705339986766\\
5049	0.099753215772101\\
9864	0.0115843048346683\\
17073	0.00962777619140532\\
27162	0.00320654317740308\\
40617	0.00134953180522152\\
57924	0.000645008353745494\\
79569	0.000344494460935705\\
106038	0.000199024565748933\\
137817	0.000122536708602598\\
175392	7.93097559679519e-05\\
219249	5.34862768927316e-05\\
269874	3.73173135919749e-05\\
327753	2.67934330699473e-05\\
393372	1.97132246950524e-05\\
467217	1.48127100311786e-05\\
};
\addplot [color=mycolor3, mark=o, mark size=2.0pt, mark options={solid, mycolor3}]
table[row sep=crcr]{
657	1.00166497368851\\
2142	0.154075925770724\\
5049	0.276540832322058\\
9864	0.0381500041118904\\
17073	0.0309599499391015\\
27162	0.011897915638109\\
40617	0.00553238781222152\\
57924	0.0028792614439328\\
79569	0.00164805580474096\\
106038	0.00101083040467783\\
137817	0.000655281300652568\\
175392	0.000443845791532779\\
219249	0.000311659573311299\\
269874	0.000225477456831667\\
327753	0.000167294212317736\\
393372	0.000126827850474004\\
467217	9.79566022291169e-05\\
};
\end{axis}
\end{tikzpicture}
\caption{Relative $\mathcal{L}^2$-error on $\Omega^{\left( j \right)}$.} \label{figure:4cubes_non_conf_2}
\end{subfigure} \\
\begin{subfigure}[t]{0.49\textwidth}
\centering
\begin{tikzpicture}[font=\footnotesize]
\begin{axis}[
width=1.0\textwidth,
height=1.0\textwidth,
xmin=-30000.0,
xmax=500000.0,
xlabel={$N_{DoFs}$},
ymin=0,
ymax=150,
ytick = { 10, 20,  30, 40, 50, 60, 70, 80, 90, 100, 110, 120, 130, 140},
ylabel={$\lvert \text{Iterations} \rvert$},
axis background/.style={fill=white}
]
\addplot [color=mycolor1, mark=*, mark size=1.0pt, mark options={solid, mycolor1}]
table[row sep=crcr]{
657	15\\
2142	12\\
5049	16\\
9864	27\\
17073	28\\
27162	36\\
40617	37\\
57924	69\\
79569	95\\
106038	89\\
137817	112\\
175392	109\\
219249	125\\
269874	86\\
327753	66\\
393372	101\\
467217	87\\
};
\addplot [color=mycolor2, mark=*, mark size=1.0pt, mark options={solid, mycolor2}]
table[row sep=crcr]{
657	15\\
2142	20\\
5049	29\\
9864	34\\
17073	34\\
27162	37\\
40617	39\\
57924	49\\
79569	55\\
106038	58\\
137817	60\\
175392	69\\
219249	72\\
269874	75\\
327753	77\\
393372	78\\
467217	79\\
};
\addplot [color=mycolor3, mark=*, mark size=1.0pt, mark options={solid, mycolor3}]
table[row sep=crcr]{
657	19\\
2142	34\\
5049	44\\
9864	49\\
17073	57\\
27162	59\\
40617	68\\
57924	72\\
79569	77\\
106038	89\\
137817	95\\
175392	100\\
219249	111\\
269874	116\\
327753	119\\
393372	129\\
467217	131\\
};
\end{axis}
\end{tikzpicture}
\caption{Number of iterations.} \label{figure:4cubes_non_conf_3}
\end{subfigure}
\hfill
\begin{subfigure}[t]{0.49\textwidth}
\centering
\begin{tikzpicture}[font=\footnotesize]

\begin{axis}[
width=1.0\textwidth,
height=1.0\textwidth,
xmin=-30000.0,
xmax=500000.0,
xlabel={$N_{DoFs}$},
ymode=log,
ymin=0,
ymax=2000,
ylabel={$s$},
yminorticks=true,
axis background/.style={fill=white}
]
\addplot [color=mycolor1, mark=*, mark size=1.0pt, mark options={solid, mycolor1}]
table[row sep=crcr]{
657	0.890036\\
2142	0.919664\\
5049	1.278806\\
9864	2.149092\\
17073	2.657286\\
27162	4.021917\\
40617	5.093085\\
57924	11.575208\\
79569	21.038837\\
106038	22.59593\\
137817	35.046511\\
175392	44.010992\\
219249	62.993494\\
269874	46.469204\\
327753	38.237397\\
393372	76.127033\\
467217	70.234128\\
};
\addplot [color=mycolor2, mark=*, mark size=1.0pt, mark options={solid, mycolor2}]
table[row sep=crcr]{
657	0.818238\\
2142	1.33025\\
5049	2.519977\\
9864	4.038817\\
17073	5.325049\\
27162	7.496528\\
40617	9.838444\\
57924	16.200243\\
79569	22.68586\\
106038	28.816142\\
137817	35.915575\\
175392	58.261211\\
219249	70.616254\\
269874	87.158946\\
327753	105.589386\\
393372	112.566106\\
467217	124.530794\\
};
\addplot [color=mycolor3, mark=*, mark size=1.0pt, mark options={solid, mycolor3}]
table[row sep=crcr]{
657	1.036623\\
2142	2.481326\\
5049	4.757451\\
9864	7.743141\\
17073	12.468757\\
27162	18.279969\\
40617	26.923099\\
57924	38.622124\\
79569	50.616595\\
106038	70.696248\\
137817	92.797727\\
175392	125.453968\\
219249	156.997041\\
269874	190.205213\\
327753	225.044054\\
393372	261.661809\\
467217	303.423997\\
};
\end{axis}
\end{tikzpicture}
\caption{Total time in seconds $s$.} \label{figure:4cubes_non_conf_4}
\end{subfigure}
\par\smallskip
\hspace{3em}
\centering
\begin{subfigure}[t]{0.9\textwidth}
\centering
\begin{tikzpicture}[font=\footnotesize]
\begin{axis}[
hide axis,
xmin=10,
xmax=50,
ymin=0,
ymax=0.4,
legend columns=3, 
legend style={draw=white!15!black,legend cell align=left, /tikz/every even column/.append style={column sep=0.4cm}}
]
\addlegendimage{mycolor1, mark=*, mark size=1.0pt, mark options={solid, mycolor1}}
\addlegendentry{$\varepsilon = 10^{-4}$};
\addlegendimage{mycolor2, mark=*, mark size=1.0pt, mark options={solid, mycolor2}}
\addlegendentry{$\varepsilon = 10^{-6}$};
\addlegendimage{mycolor3, forget plot}
\addlegendentry{$\varepsilon = 10^{-8}$};
\addlegendimage{mycolor3, mark=triangle, mark size=2.0pt, mark options={solid, mycolor3}}
\addlegendentry{patch 1};
\addlegendimage{mycolor3, mark=+, mark size=2.0pt, mark options={solid, mycolor3}}
\addlegendentry{patch 2};
\addlegendimage{mycolor3, mark=o, mark size=2.0pt, mark options={solid, mycolor3}}
\addlegendentry{patch 3};
\end{axis}
\end{tikzpicture}
\end{subfigure}
\caption{Performance for different refinements and tolerances on the nonconforming multi-patch geometry with $4$ cuboids shown in~\ref{figure:4cubes_conf_1} with~\eqref{eq:4cubes_sol} as analytical solution.} \label{figure:4cubes_non_conf_all} 
\end{figure}

\begin{figure}
\definecolor{mycolor1}{rgb}{0.00000,0.44700,0.74100}
\definecolor{mycolor2}{rgb}{0.85000,0.32500,0.09800}
\definecolor{mycolor3}{rgb}{0.92900,0.69400,0.12500}
\definecolor{mycolor4}{rgb}{0.49400,0.18400,0.55600}
\definecolor{mycolor5}{rgb}{0.46600,0.67400,0.18800}
\begin{subfigure}[t]{0.49\textwidth}
\centering
\includegraphics[width=0.80\textwidth]{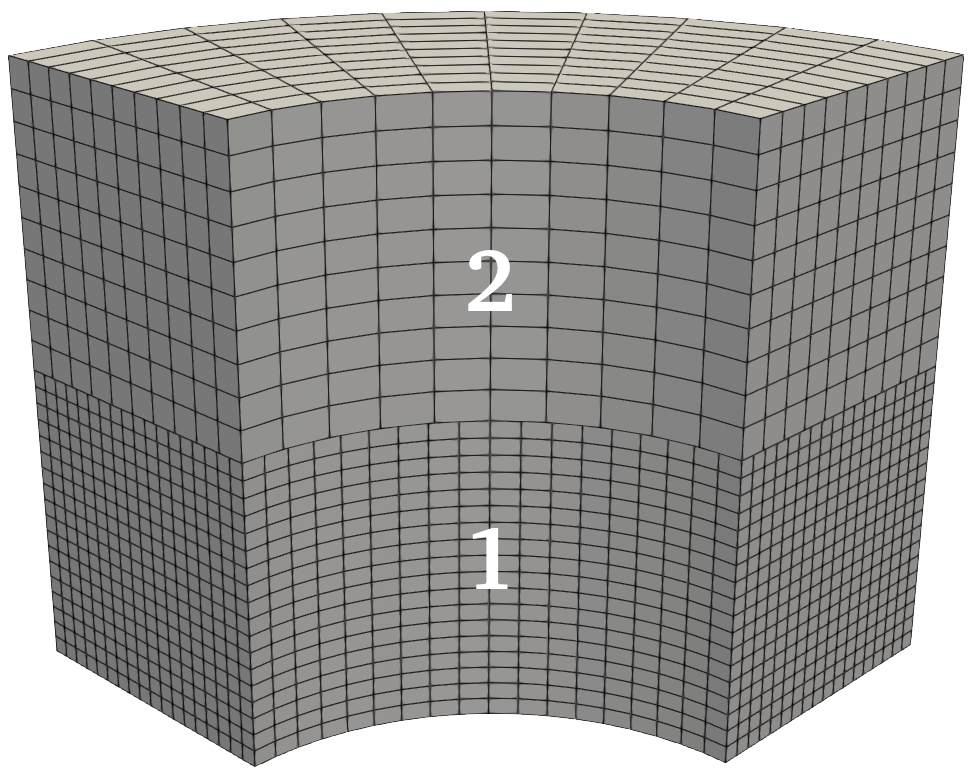}
\caption{2 annuli placed on top of each other.}
\label{figure:2annulus_conf_1}
\end{subfigure}
\hfill
\begin{subfigure}[t]{0.49\textwidth}
\centering
\begin{tikzpicture}[font=\footnotesize]
\begin{axis}[
width=1.0\textwidth,
height=1.0\textwidth,
xmin=-50000,
xmax=950000,
xlabel={$N_{DoFs}$},
ymode=log,
ymin=1e-09,
ymax=10,
yminorticks=true,
yminorgrids=true,
ymajorgrids=true,
xminorgrids=true,
ytick = {1e-4, 1e-6, 1e-8},
extra y ticks = {1e-2, 1e-3, 1e-5, 1e-7},
extra y tick style={grid=none},
ylabel={$R ( y_{h}, y_{sol} )$},
axis background/.style={fill=white}
]
\addplot [color=mycolor1, mark=*, mark size=1.0pt, mark options={solid, mycolor1}]
table[row sep=crcr]{
128	0.0614838021605079\\
686	0.000714919509814133\\
2000	7.61714357351822e-05\\
4394	1.85778491700005e-05\\
8192	5.32880820465821e-05\\
13718	1.96012379631342e-05\\
21296	2.03073353552473e-05\\
31250	6.61834531874783e-05\\
43904	0.000115427828527827\\
59582	0.000115168556212762\\
78608	0.000103472162417666\\
101306	0.000132436022163884\\
128000	0.00015743496578676\\
159014	0.000279354868902368\\
194672	0.000181032643223315\\
235298	0.00023625846970616\\
281216	0.000246135905438495\\
332750	0.000201475793776896\\
390224	0.000183895074046186\\
453962	0.000318486076223687\\
524288	0.000362773309938068\\
601526	0.000279167323954689\\
686000	0.000387254305238004\\
778034	0.000389158384155277\\
877952	0.000282138583809976\\
};
\addplot [color=mycolor2, mark=*, mark size=1.0pt, mark options={solid, mycolor2}]
table[row sep=crcr]{
128	0.0614834083438276\\
686	0.000714972562910947\\
2000	7.61758690579522e-05\\
4394	1.85178356009441e-05\\
8192	6.55308604082242e-06\\
13718	2.87949904782329e-06\\
21296	1.49660506241484e-06\\
31250	9.34295166379291e-07\\
43904	2.92526547815184e-06\\
59582	1.71608694889159e-06\\
78608	1.7310643060766e-06\\
101306	2.14288811849946e-06\\
128000	3.37769658233796e-06\\
159014	1.37948513119367e-06\\
194672	3.05411823408374e-06\\
235298	3.29071428261643e-06\\
281216	4.50134422472026e-06\\
332750	4.34107986398998e-06\\
390224	2.75579156374681e-06\\
453962	1.53161983051861e-06\\
524288	2.56376103428688e-06\\
601526	2.54035352741086e-06\\
686000	2.56827305197312e-06\\
778034	4.61315196196826e-06\\
877952	2.89171740015336e-06\\
};
\addplot [color=mycolor3, mark=*, mark size=1.0pt, mark options={solid, mycolor3}]
table[row sep=crcr]{
128	0.061483409436566\\
686	0.000714976870054499\\
2000	7.61758729332798e-05\\
4394	1.85178512522541e-05\\
8192	6.5456728054634e-06\\
13718	2.87169762226226e-06\\
21296	1.45112843856462e-06\\
31250	8.10241189944509e-07\\
43904	4.87303774194017e-07\\
59582	3.10497193989661e-07\\
78608	2.0786155548068e-07\\
101306	1.44033473908481e-07\\
128000	1.07589831234533e-07\\
159014	7.99753984814963e-08\\
194672	7.10878198765579e-08\\
235298	4.55413165656174e-08\\
281216	3.92381918678057e-08\\
332750	5.03871248403483e-08\\
390224	4.06467194109856e-08\\
453962	3.75467117543191e-08\\
524288	4.26804805802356e-08\\
601526	2.95371103373473e-08\\
686000	3.08942190048982e-08\\
778034	2.94596550238291e-08\\
877952	4.34832052459791e-08\\
};
\addplot [color=mycolor5, mark=star, mark size=2.0pt, mark options={solid, mycolor5}]
table[row sep=crcr]{
128	0.0615027999746377\\
686	0.000714976759039962\\
2000	7.61758706251123e-05\\
4394	1.85178578709077e-05\\
8192	6.54567672883507e-06\\
13718	2.87169314933354e-06\\
21296	1.45112091840995e-06\\
31250	8.10235733898558e-07\\
43904	4.87277236122385e-07\\
59582	3.10378099458291e-07\\
78608	2.06958472743454e-07\\
101306	1.43247079360756e-07\\
};
\end{axis}
\end{tikzpicture}
\caption{Relative $\mathcal{L}^2$-error on $\Omega$.}\label{figure:2annulus_conf_2}
\end{subfigure} \\
\begin{subfigure}[t]{0.49\textwidth}
\centering
\begin{tikzpicture}[font=\footnotesize]

\begin{axis}[
width=1.0\textwidth,
height=1.0\textwidth,
xmin=-50000,
xmax=950000,
xlabel={$N_{DoFs}$},
ymin=0,
ymax=25,
ylabel={$\lvert \text{Iterations} \rvert$},
ytick={0, 5, 10, 15, 20},
axis background/.style={fill=white}
]
\addplot [color=mycolor1, mark=*, mark size=1.0pt, mark options={solid, mycolor1}]
table[row sep=crcr]{
128	2\\
686	3\\
2000	5\\
4394	4\\
8192	6\\
13718	5\\
21296	4\\
31250	4\\
43904	5\\
59582	3\\
78608	4\\
101306	4\\
128000	4\\
159014	4\\
194672	4\\
235298	2\\
281216	3\\
332750	3\\
390224	5\\
453962	3\\
524288	4\\
601526	4\\
686000	4\\
778034	4\\
877952	4\\
};
\addplot [color=mycolor2, mark=*, mark size=1.0pt, mark options={solid, mycolor2}]
table[row sep=crcr]{
128	3\\
686	6\\
2000	8\\
4394	9\\
8192	10\\
13718	10\\
21296	10\\
31250	10\\
43904	8\\
59582	8\\
78608	10\\
101306	10\\
128000	9\\
159014	10\\
194672	11\\
235298	10\\
281216	8\\
332750	9\\
390224	12\\
453962	10\\
524288	12\\
601526	12\\
686000	9\\
778034	10\\
877952	11\\
};
\addplot [color=mycolor3, mark=*, mark size=1.0pt, mark options={solid, mycolor3}]
table[row sep=crcr]{
128	4\\
686	7\\
2000	11\\
4394	12\\
8192	13\\
13718	15\\
21296	15\\
31250	16\\
43904	16\\
59582	15\\
78608	16\\
101306	16\\
128000	15\\
159014	17\\
194672	18\\
235298	16\\
281216	18\\
332750	16\\
390224	20\\
453962	20\\
524288	18\\
601526	20\\
686000	22\\
778034	18\\
877952	22\\
};
\end{axis}
\end{tikzpicture}
\caption{Number of iterations.}\label{figure:2annulus_conf_3}
\end{subfigure}
\hfill
\begin{subfigure}[t]{0.49\textwidth}
\centering
\begin{tikzpicture}[font=\footnotesize]

\begin{axis}[
width=1.0\textwidth,
height=1.0\textwidth,
xmin=-50000,
xmax=950000,
xlabel={$N_{DoFs}$},
ymode=log,
ymin=0,
ymax=2000,
yminorticks=true,
ylabel={$s$},
axis background/.style={fill=white}
]
\addplot [color=mycolor1, mark=*, mark size=1.0pt, mark options={solid, mycolor1}]
table[row sep=crcr]{
128	0.922179\\
686	0.472032\\
2000	0.47249\\
4394	0.539618\\
8192	0.802404\\
13718	1.040031\\
21296	1.4964\\
31250	2.146399\\
43904	3.341468\\
59582	4.801441\\
78608	6.029267\\
101306	8.088582\\
128000	10.429159\\
159014	13.014364\\
194672	15.784305\\
235298	20.43132\\
281216	24.756368\\
332750	30.101685\\
390224	36.500118\\
453962	45.355325\\
524288	51.875999\\
601526	62.619312\\
686000	70.817721\\
778034	85.36523\\
877952	99.927727\\
};
\addplot [color=mycolor2, mark=*, mark size=1.0pt, mark options={solid, mycolor2}]
table[row sep=crcr]{
128	0.190962\\
686	0.339863\\
2000	0.533342\\
4394	0.712331\\
8192	1.072404\\
13718	1.560608\\
21296	2.254287\\
31250	3.298842\\
43904	4.238994\\
59582	5.842601\\
78608	7.553753\\
101306	9.505782\\
128000	11.985629\\
159014	15.09187\\
194672	18.254746\\
235298	22.953385\\
281216	27.173082\\
332750	32.699355\\
390224	40.018888\\
453962	48.478019\\
524288	56.523132\\
601526	67.382394\\
686000	73.508959\\
778034	89.141766\\
877952	101.730845\\
};
\addplot [color=mycolor3, mark=*, mark size=1.0pt, mark options={solid, mycolor3}]
table[row sep=crcr]{
128	0.202708\\
686	0.388061\\
2000	0.678638\\
4394	0.983306\\
8192	1.473624\\
13718	2.682279\\
21296	3.707958\\
31250	5.773554\\
43904	7.813065\\
59582	9.257283\\
78608	11.780984\\
101306	14.499138\\
128000	17.09105\\
159014	22.425681\\
194672	25.9576\\
235298	30.649424\\
281216	36.360086\\
332750	41.027334\\
390224	53.933574\\
453962	63.269286\\
524288	68.887082\\
601526	82.595157\\
686000	95.950952\\
778034	105.222687\\
877952	125.770819\\
};
\addplot [color=mycolor5, mark=star, mark size=2.0pt, mark options={solid, mycolor5}]
table[row sep=crcr]{
128	0.772802\\
686	0.750849\\
2000	2.760513\\
4394	7.167908\\
8192	15.457038\\
13718	30.119744\\
21296	52.356997\\
31250	84.98963\\
43904	138.15897\\
59582	209.897785\\
78608	303.351681\\
101306	513.304024\\
};
\end{axis}

\end{tikzpicture}
\caption{Total time in seconds $s$.} \label{figure:2annulus_conf_4}
\end{subfigure}
\par \smallskip
\hspace{3em}
\begin{subfigure}[t]{0.9\textwidth}
\centering
\begin{tikzpicture}[font=\footnotesize]
\begin{axis}[
hide axis,
xmin=10,
xmax=50,
ymin=0,
ymax=0.4,
legend columns=5, 
legend style={draw=white!15!black,legend cell align=left, /tikz/every even column/.append style={column sep=0.4cm}}
]
\addlegendimage{mycolor1, mark=*, mark size=1.0pt, mark options={solid, mycolor1}}
\addlegendentry{$\varepsilon = 10^{-4}$};
\addlegendimage{mycolor2, mark=*, mark size=1.0pt, mark options={solid, mycolor2}}
\addlegendentry{$\varepsilon = 10^{-6}$};
\addlegendimage{mycolor3, mark=*, mark size=1.0pt, mark options={solid, mycolor3}}
\addlegendentry{$\varepsilon = 10^{-8}$};
\addlegendimage{color=mycolor5, mark=star, mark size=2.0pt, mark options={solid, mycolor5}}
\addlegendentry{\textsc{GeoPDEs 3.0}};
\end{axis}
\end{tikzpicture}
\end{subfigure}
\caption{Performance for different refinements and tolerances on the conforming multi-patch geometry with $2$ annuli shown in~\ref{figure:2annulus_conf_1} with~\eqref{eq:2annulus_sol} as analytical solution.} \label{figure:2annulus_conf_all} 
\end{figure}

The last two test cases are studied on a NURBS geometry shown in Figure \ref{figure:2annulus_conf_1} consisting of two vertically stacked annuli. The analytical is given by:
\begin{equation}
\label{eq:2annulus_sol}
y_{sol} \left( x, y, z \right) = \left( x^2 + y^2 - 1 \right) \, \left( x^2 + y^2 - 4 \right) \, x \, y \, z \, \left( z - 2 \right).
\end{equation}
We use B-splines of degree $3$ in both the conforming case (results in Figure \ref{figure:2annulus_conf_all}) and the nonconforming case (results in Figure \ref{figure:2annulus_non_conf_all}). In the conforming case, we start with $N^{\left( j \right)}_{DoFs} = 64$ and end with $N^{\left( j \right)}_{DoFs} = 438.976$ for $j = 1, 2$. In the nonconforming case, we start with $N^{\left( 1 \right)}_{DoFs} = 125$, $N^{\left( 2 \right)}_{DoFs} = 64$ and end with $N^{\left( 1 \right)}_{DoFs} = 456.533$, $N^{\left( 2 \right)}_{DoFs}= 64.000$. We observe similar results here too.

\begin{figure}
\definecolor{mycolor1}{rgb}{0.00000,0.44700,0.74100}
\definecolor{mycolor2}{rgb}{0.85000,0.32500,0.09800}
\definecolor{mycolor3}{rgb}{0.92900,0.69400,0.12500}
\definecolor{mycolor4}{rgb}{0.49400,0.18400,0.55600}
\definecolor{mycolor5}{rgb}{0.46600,0.67400,0.18800}

\begin{subfigure}[t]{0.49\textwidth}
\centering
\begin{tikzpicture}[font=\footnotesize]
\begin{axis}[
width=1.0\textwidth,
height=1.0\textwidth,
xmin=-50000.0,
xmax=550000.0,
xlabel={$N_{DoFs}$},
ymode=log,
ylabel={$R ( y_{h}, y_{sol} )$},
ymin=1e-09,
ymax=10,
yminorticks=true,
yminorgrids=true,
ymajorgrids=true,
xminorgrids=true,
ytick = {1e-4, 1e-6, 1e-8},
extra y ticks = {1e-2, 1e-3, 1e-5, 1e-7},
extra y tick style={grid=none},
axis background/.style={fill=white}
]
\addplot [color=mycolor1, mark=*, mark size=1.0pt, mark options={solid, mycolor1}]
table[row sep=crcr]{
189	0.0357306902069579\\
1674	0.000188585308622369\\
5913	2.74584729703464e-05\\
14364	6.49994200330134e-05\\
28485	7.73663173726182e-05\\
49734	8.35921670814048e-05\\
79569	0.000119726668138224\\
119448	0.000100020251342103\\
170829	0.000167656422429396\\
235170	0.000149131186983833\\
313929	0.000199290857801499\\
408564	0.000213147497198728\\
520533	0.000314722287671984\\
};
\addplot [color=mycolor2, mark=*, mark size=1.0pt, mark options={solid, mycolor2}]
table[row sep=crcr]{
189	0.035708644673771\\
1674	0.000289267779322765\\
5913	2.4467117960616e-05\\
14364	5.46609269220151e-06\\
28485	2.08738555326337e-06\\
49734	1.51725049647835e-06\\
79569	3.18887766719666e-06\\
119448	3.465403088529e-06\\
170829	2.84173411581339e-06\\
235170	1.94594418631383e-06\\
313929	2.82619052161857e-06\\
408564	3.11063074482129e-06\\
520533	5.00117628844007e-06\\
};
\addplot [color=mycolor3, mark=*, mark size=1.0pt, mark options={solid, mycolor3}]
table[row sep=crcr]{
189	0.0357087845134815\\
1674	0.000289355238728026\\
5913	2.44970436547983e-05\\
14364	5.25629670552451e-06\\
28485	1.7137874765453e-06\\
49734	7.09626147641871e-07\\
79569	3.45694322497043e-07\\
119448	1.88081613893134e-07\\
170829	1.13367693145323e-07\\
235170	8.93367416370995e-08\\
313929	5.70453997411176e-08\\
408564	3.75387101805494e-08\\
520533	4.93418389081636e-08\\
};
\end{axis}
\end{tikzpicture}
\caption{Relative $\mathcal{L}^2$-error on $\Omega$.} \label{figure:2annulus_non_conf_2}
\end{subfigure}
\hfill
\begin{subfigure}[t]{0.49\textwidth}
\centering
\begin{tikzpicture}[font=\footnotesize]
\begin{axis}[
width=1.0\textwidth,
height=1.0\textwidth,
xmin=-50000.0,
xmax=550000.0,
xlabel={$N_{DoFs}$},
ymode=log,
ymin=1.0e-09,
ymax=10,
ylabel={$R^{ (j) } ( y_h, y_{sol} )$},
yminorticks=true,
yminorgrids=true,
ymajorgrids=true,
xminorgrids=true,
ytick = {1e-8},
extra y ticks = {1e-2, 1e-4, 1e-3, 1e-5, 1e-6, 1e-7},
extra y tick style={grid=none},
axis background/.style={fill=white},
]
\addplot [color=mycolor3, mark size=2.0pt, mark=+, mark options={solid, mycolor3}]
table[row sep=crcr]{
189	0.0226578595403781\\
1674	0.000179841354813607\\
5913	1.39797212842101e-05\\
14364	2.86177619113239e-06\\
28485	9.02316421452351e-07\\
49734	3.63751332504718e-07\\
79569	1.74904705527729e-07\\
119448	9.4424304507564e-08\\
170829	5.82205024191012e-08\\
235170	5.72499654245143e-08\\
313929	3.54407028427331e-08\\
408564	2.19072745224815e-08\\
520533	4.17367016294841e-08\\
};
\addplot [color=mycolor3, mark size=2.0pt, mark=o, mark options={solid, mycolor3}]
table[row sep=crcr]{
189	0.0611988723515742\\
1674	0.000714320193509632\\
5913	7.61686484614981e-05\\
14364	1.85171665728931e-05\\
28485	6.5455666666483e-06\\
49734	2.8716589073744e-06\\
79569	1.45115931039544e-06\\
119448	8.10403715162767e-07\\
170829	4.87372262831355e-07\\
235170	3.10542810978999e-07\\
313929	2.08001788040498e-07\\
408564	1.47989478502558e-07\\
520533	1.03591778552654e-07\\
};
\end{axis}
\end{tikzpicture}
\caption{Relative $\mathcal{L}^2$-error on $\Omega^{\left( j \right)}$.} \label{figure:2annulus_non_conf_1}
\end{subfigure} \\
\begin{subfigure}[t]{0.49\textwidth}
\centering
\begin{tikzpicture}[font=\footnotesize]
\begin{axis}[
width=1.0\textwidth,
height=1.0\textwidth,
xmin=-50000.0,
xmax=550000.0,
xlabel={$N_{DoFs}$},
ymin=0,
ymax=55,
ytick = {10, 20, 30, 40, 50},
ylabel={$\lvert \text{Iterations} \rvert$},
axis background/.style={fill=white}
]
\addplot [color=mycolor1, mark=*, mark size=1.0pt, mark options={solid, mycolor1}]
table[row sep=crcr]{
189	3\\
1674	6\\
5913	8\\
14364	9\\
28485	8\\
49734	4\\
79569	7\\
119448	2\\
170829	4\\
235170	4\\
313929	3\\
408564	3\\
520533	4\\
};
\addplot [color=mycolor2, mark=*, mark size=1.0pt, mark options={solid, mycolor2}]
table[row sep=crcr]{
189	4\\
1674	10\\
5913	13\\
14364	16\\
28485	18\\
49734	17\\
79569	19\\
119448	18\\
170829	27\\
235170	19\\
313929	18\\
408564	16\\
520533	19\\
};
\addplot [color=mycolor3, mark=*, mark size=1.0pt, mark options={solid, mycolor3}]
table[row sep=crcr]{
189	5\\
1674	14\\
5913	18\\
14364	24\\
28485	27\\
49734	29\\
79569	33\\
119448	35\\
170829	36\\
235170	38\\
313929	40\\
408564	46\\
520533	48\\
};
\end{axis}
\end{tikzpicture}
\caption{Number of iterations.} \label{figure:2annulus_non_conf_3}
\end{subfigure}
\hfill
\begin{subfigure}[t]{0.49\textwidth}
\centering
\begin{tikzpicture}[font=\footnotesize]
\begin{axis}[
width=1.0\textwidth,
height=1.0\textwidth,
xmin=-50000.0,
xmax=550000.0,
xlabel={$N_{DoFs}$},
ymode=log,
ymin=0,
ymax=500,
ylabel={$s$},
yminorticks=true,
axis background/.style={fill=white}
]
\addplot [color=mycolor1, mark=*, mark size=1.0pt, mark options={solid, mycolor1}]
table[row sep=crcr]{
189	0.255205\\
1674	0.450065\\
5913	0.611787\\
14364	1.158111\\
28485	2.226402\\
49734	3.59633\\
79569	6.091714\\
119448	9.020547\\
170829	13.707503\\
235170	19.912482\\
313929	27.864461\\
408564	37.353422\\
520533	52.504072\\
};
\addplot [color=mycolor2, mark=*, mark size=1.0pt, mark options={solid, mycolor2}]
table[row sep=crcr]{
189	0.233583\\
1674	0.450605\\
5913	0.895542\\
14364	1.812433\\
28485	3.85491\\
49734	5.757985\\
79569	9.731797\\
119448	13.001196\\
170829	22.451869\\
235170	25.545399\\
313929	33.808657\\
408564	42.34857\\
520533	59.281693\\
};
\addplot [color=mycolor3, mark=*, mark size=1.0pt, mark options={solid, mycolor3}]
table[row sep=crcr]{
189	0.27682\\
1674	0.657416\\
5913	1.363327\\
14364	3.066223\\
28485	6.328221\\
49734	10.746778\\
79569	17.506944\\
119448	25.396782\\
170829	34.763294\\
235170	44.67855\\
313929	62.360559\\
408564	81.335004\\
520533	105.227903\\
};
\end{axis}

\end{tikzpicture}
\caption{Total time in seconds $s$.} \label{figure:2annulus_non_conf_4}
\end{subfigure}
\par\smallskip
\hspace{3em}
\centering
\begin{subfigure}[t]{0.9\textwidth}
\centering
\begin{tikzpicture}[font=\footnotesize]
\begin{axis}[
hide axis,
xmin=10,
xmax=50,
ymin=0,
ymax=0.4,
legend columns=3, 
legend style={draw=white!15!black,legend cell align=left, /tikz/every even column/.append style={column sep=0.4cm}}
]
\addlegendimage{mycolor1, mark=*, mark size=1.0pt, mark options={solid, mycolor1}}
\addlegendentry{$\varepsilon = 10^{-4}$};
\addlegendimage{mycolor2, mark=*, mark size=1.0pt, mark options={solid, mycolor2}}
\addlegendentry{$\varepsilon = 10^{-6}$};
\addlegendimage{mycolor3, forget plot}
\addlegendentry{$\varepsilon = 10^{-8}$};
\addlegendimage{mycolor3, mark=+, mark size=2.0pt, mark options={solid, mycolor3}}
\addlegendentry{patch 1};
\addlegendimage{mycolor3, mark=o, mark size=2.0pt, mark options={solid, mycolor3}}
\addlegendentry{patch 2};
\end{axis}
\end{tikzpicture}
\end{subfigure}
\caption{Performance for different refinements and tolerances on the nonconforming multi-patch geometry with $2$ annuli shown in~\ref{figure:2annulus_conf_1} with~\eqref{eq:2annulus_sol} as analytical solution.} \label{figure:2annulus_non_conf_all} 
\end{figure}

\subsection{Optimal control}

We now illustrate the performance of our low-rank approach in section~\ref{section:optimization} for the optimal control problem~\eqref{eq:optimization1} to~\eqref{eq:optimization4} using the following desired state
\begin{equation}
\label{eq:OC_sol}
\hat{y}\left( x, y, z; t \right) = e^{\frac{1}{t+1}} \, y_{sol} \left( x, y, z \right).
\end{equation}
and we set $T = 1$ and $N_t = 10$. Our main goal is to illustrate the robustness of our method with respect to changes in the penalty parameter $\alpha$. We solve~\eqref{eq:OC_Schur_complement} using \texttt{tt\_gmres\_block.m}, where we set \texttt{max\_iters} to $10$, \texttt{restart} to $20$, and \texttt{tol} to $10^{-5}$. We solve the inner linear systems $\bar{\mathcal{K}}$ and $\bar{\mathcal{K}}^{\top}$ with the approach described in section~\ref{section:IETI_method}, i.e.~for solving~\eqref{eq:preconditioned_Schu_complement} we use \texttt{tt\_gmres\_block.m}\footnote{with parameters \texttt{max\_iters} $ = 10$, \texttt{restart} $= 20$ and \texttt{tol} $= 10^{-6}$} and for solving the local linear systems inside~\eqref{eq:preconditioned_Schu_complement} we use \texttt{amen\_solve2.m}\footnote{with parameters \texttt{nswp} $ = 20$, \texttt{kickrank} $= 2$ and \texttt{tol} $= 10^{-7}$}. As before we use \texttt{amen\_block\_solve.m} for solving the weight function interpolation system~\eqref{eq:weightfunction_linear_system} with a tolerance of $10^{-8}$. \\

\begin{figure}
\definecolor{mycolor1}{rgb}{0.00000,0.44700,0.74100}
\definecolor{mycolor2}{rgb}{0.85000,0.32500,0.09800}
\definecolor{mycolor3}{rgb}{0.92900,0.69400,0.12500}
\definecolor{mycolor4}{rgb}{0.49400,0.18400,0.55600}
\definecolor{mycolor5}{rgb}{0.46600,0.67400,0.18800}

\begin{subfigure}[t]{0.49\textwidth}
\centering
\begin{tikzpicture}[font=\footnotesize]
\begin{axis}[
width=0.95\textwidth,
height=1.0\textwidth,
xmin=0.5,
xmax=5.5,
ymode=log,
ymin=1.3,
ymax=3.5,
yminorticks=true,
ylabel={$\mathcal{L} \left( y, u, \mu \right)$},
xtick={1, 2, 3, 4, 5},
xticklabels={1e-4, 1e-3, 1e-2, 1e-1, 1},
xlabel={$\alpha$},
axis background/.style={fill=white}
]
\addplot [mycolor1, mark=o, mark size=2.0pt, mark options={solid, mycolor1}]
table[row sep=crcr]{
1	1.493832674190639\\
2	2.505335763703290\\
3	2.687617082033170\\
4	2.707320290915187\\
5	2.709306571914217\\
};
\addplot [mycolor2, mark=*, mark size=1.0pt, mark options={solid, mycolor2}]
table[row sep=crcr]{
1	1.49383299909775\\
2	2.50533585048723\\
3	2.68761708630727\\
4	2.70732028553885\\
5	2.70930656555751\\
};
\end{axis}
\end{tikzpicture}
\caption{Value of objective function $\mathcal{L} \left( y, u, \mu \right)$.} \label{figure:OC_4cubes_1}
\end{subfigure}
\hfill
\begin{subfigure}[t]{0.49\textwidth}
\centering
\begin{tikzpicture}[font=\footnotesize]
\begin{axis}[
width=0.95\textwidth,
height=1.0\textwidth,
xmin=0.5,
xmax=5.5,
ymode=log,
ymin=0,
ymax=1e6,
yminorticks=true,
ylabel={$\lVert u \rVert$},
xtick={1, 2, 3, 4, 5},
xticklabels={1e-4, 1e-3, 1e-2, 1e-1, 1},
xlabel={$\alpha$},
axis background/.style={fill=white}
]
\addplot [mycolor1, mark=o, mark size=2.0pt, mark options={solid, mycolor1}]
table[row sep=crcr]{
1	9.820308198156809e+04\\
2	1.648856394448072e+04\\
3	1.769254987716603e+03\\
4	1.782273805697482e+02\\
5	17.835862920892509\\
};
\addplot [mycolor2, mark=*, mark size=1.0pt, mark options={solid, mycolor2}]
table[row sep=crcr]{
1	5.479898787436817e+04\\
2	9.200898035957045e+03\\
3	9.872742162947811e+02\\
4	99.453893227634651\\
5	9.952713070615829\\
};
\end{axis}
\end{tikzpicture}
\caption{Norm of control $\lVert u \rVert$.} \label{figure:OC_4cubes_2}
\end{subfigure} \\
\begin{subfigure}[t]{0.49\textwidth}
\centering
\begin{tikzpicture}[font=\footnotesize]
\begin{axis}[
width=0.95\textwidth,
height=1.0\textwidth,
xmin=0.5,
xmax=5.5,
ymin=1,
ymax=5,
ylabel={$\lvert \text{Iterations} \rvert$},
xtick={1, 2, 3, 4, 5},
xticklabels={1e-4, 1e-3, 1e-2, 1e-1, 1},
xlabel={$\alpha$},
axis background/.style={fill=white}
]
\addplot [mycolor1, mark=o, mark size=2.0pt, mark options={solid, mycolor1}]
table[row sep=crcr]{
1	4\\
2	3\\
3	2\\
4	2\\
5	2\\
};
\addplot [mycolor2, mark=*, mark size=1.0pt, mark options={solid, mycolor2}]
table[row sep=crcr]{
1	4\\
2	3\\
3	2\\
4	2\\
5	2\\
};
\end{axis}
\end{tikzpicture}
\caption{Number of iterations.} \label{figure:OC_4cubes_3}
\end{subfigure}
\hfill
\begin{subfigure}[t]{0.49\textwidth}
\centering
\begin{tikzpicture}[font=\footnotesize]
\begin{axis}[
width=0.95\textwidth,
height=1.0\textwidth,
xmin=0.5,
xmax=5.5,
ymode=log,
ymin=5500,
ymax=1.3e+04,
yminorticks=true,
ylabel={$s$},
xtick={1, 2, 3, 4, 5},
xticklabels={1e-4, 1e-3, 1e-2, 1e-1, 1},
xlabel={$\alpha$},
axis background/.style={fill=white}
]
\addplot [mycolor1, mark=o, mark size=2.0pt, mark options={solid, mycolor1}]
table[row sep=crcr]{
1	1.178120855300000e+04\\
2	7.351027785000000e+03\\
3	6.945888416000000e+03\\
4	6.226359552000000e+03\\
5	6.500587958000000e+03\\
};
\addplot [mycolor2, mark=*, mark size=1.0pt, mark options={solid, mycolor2}]
table[row sep=crcr]{
1	1.028154593300000e+04\\
2	8.529266267999999e+03\\
3	6.775572543000000e+03\\
4	6.789527227000000e+03\\
5	6.669345227000000e+03\\
};
\end{axis}
\end{tikzpicture}
\caption{Total time in seconds $s$.} \label{figure:OC_4cubes_4}
\end{subfigure}
\par\smallskip
\hspace{3em}
\centering
\begin{subfigure}[t]{0.9\textwidth}
\centering
\begin{tikzpicture}[font=\footnotesize]
\begin{axis}[
hide axis,
xmin=10,
xmax=50,
ymin=0,
ymax=0.4,
legend columns=2, 
legend style={draw=white!15!black,legend cell align=left, /tikz/every even column/.append style={column sep=0.4cm}}
]
\addlegendimage{mycolor1, mark=o, mark size=2.0pt, mark options={solid, mycolor1}}
\addlegendentry{conforming};
\addlegendimage{mycolor2, mark=*, mark size=1.0pt, mark options={solid, mycolor2}}
\addlegendentry{nonconforming};
\end{axis}
\end{tikzpicture}
\end{subfigure}
\caption{Stability of the method presented in section~\ref{section:optimization} depending on the parameter $\alpha$ on the multi-patch geometry with $4$ cuboids shown in~\ref{figure:4cubes_conf_1} for the conforming and not nonconforming case.} \label{figure:OC_4cubes_all} 
\end{figure}
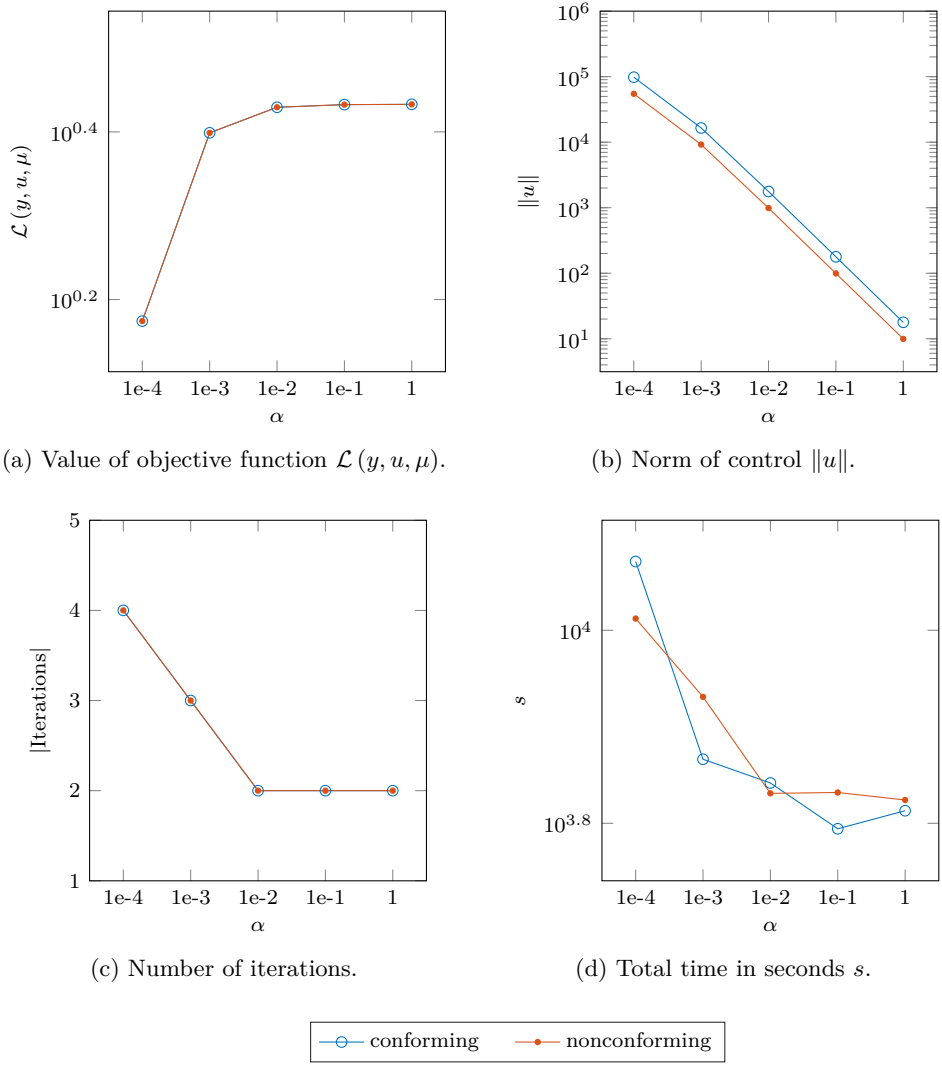

We show in Figure \ref{figure:OC_4cubes_all} the results for the multi-patch geometry shown in Figure \ref{figure:4cubes_conf_1} with desired state defined by~\eqref{eq:OC_sol} where $y_{sol} \left( x, y, z \right)$ is defined by~\eqref{eq:4cubes_sol}. We use B-splines of degree $3$ in both the conforming case and the nonconforming case. In the conforming case we have a discretization of $N^{\left( j \right)}_{DoFs} = 373.248$ DoFs for $j = 1, \ldots, 4$, so that we consider a total of $14.929.920$ many DoFs for $N_t = 10$ time steps. In the nonconforming case we have the division $N^{\left( 1 \right)}_{DoFs} = 357.911$, $N^{\left( 2 \right)}_{DoFs}, N^{\left( 4 \right)}_{DoFs} = 50.653$, $N^{\left( 3 \right)}_{DoFs} = 8000$, so that we consider a total of $4.672.170$ DoFs with $N_t = 10$. We observe in both cases that the objective function value decreases while the norm of the control increases for decreasing $\alpha$. The number of iterations is small which suggests that the preconditioner approximates the full operator sufficiently and hence reducing the number of \ttgmres\ iterations. The computational times remain moderate in both cases and do only slightly vary in a similar way as the number of iterations, which is to be expected. \\

\begin{figure}
\definecolor{mycolor1}{rgb}{0.00000,0.44700,0.74100}
\definecolor{mycolor2}{rgb}{0.85000,0.32500,0.09800}
\definecolor{mycolor3}{rgb}{0.92900,0.69400,0.12500}
\definecolor{mycolor4}{rgb}{0.49400,0.18400,0.55600}
\definecolor{mycolor5}{rgb}{0.46600,0.67400,0.18800}

\begin{subfigure}[t]{0.49\textwidth}
\centering
\begin{tikzpicture}[font=\footnotesize]
\begin{axis}[
width=0.95\textwidth,
height=1.0\textwidth,
xmin=0.5,
xmax=5.5,
ymode=log,
ymin=0,
ymax=200,
yminorticks=true,
ylabel={$\mathcal{L} \left( y, u, \mu \right)$},
xtick={1, 2, 3, 4, 5},
xticklabels={1e-4, 1e-3, 1e-2, 1e-1, 1},
xlabel={$\alpha$},
axis background/.style={fill=white}
]
\addplot [mycolor1, mark=o, mark size=2.0pt, mark options={solid, mycolor1}]
table[row sep=crcr]{
1	6.96152844973917\\
2	35.7600912281175\\
3	91.8081329236656\\
4	116.894293648742\\
5	120.417224784926\\
};
\addplot [mycolor2, mark=*, mark size=1.0pt, mark options={solid, mycolor2}]
table[row sep=crcr]{
1	6.961528476612264\\
2	35.760091219048896\\
3	91.808132782057740\\
4	1.168942936655431e2\\
5	1.204172247942133e2\\
};
\end{axis}
\end{tikzpicture}
\caption{Value of objective function $\mathcal{L} \left( y, u, \mu \right)$.} \label{figure:OC_2annulus_1}
\end{subfigure}
\hfill
\begin{subfigure}[t]{0.49\textwidth}
\centering
\begin{tikzpicture}[font=\footnotesize]
\begin{axis}[
width=0.95\textwidth,
height=1.0\textwidth,
xmin=0.5,
xmax=5.5,
ymode=log,
ymin=0,
ymax=1e5,
yminorticks=true,
ylabel={$\lVert u \rVert$},
xtick={1, 2, 3, 4, 5},
xticklabels={1e-4, 1e-3, 1e-2, 1e-1, 1},
xlabel={$\alpha$},
axis background/.style={fill=white}
]
\addplot [mycolor1, mark=o, mark size=2.0pt, mark options={solid, mycolor1}]
table[row sep=crcr]{
1	64175.9600285374\\
2	39712.4908447363\\
3	12897.6785564197\\
4	1750.57739753495\\
5	181.724094773199\\
};
\addplot [mycolor2, mark=*, mark size=1.0pt, mark options={solid, mycolor2}]
table[row sep=crcr]{
1	6.962252337808002e4\\
2	4.307786577641687e4\\
3	1.399092374425819e4\\
4	1.898983365247624e3\\
5	1.971301093659045e2\\
};
\end{axis}
\end{tikzpicture}
\caption{Norm of control $\lVert u \rVert$.} \label{figure:OC_2annulus_2}
\end{subfigure} \\
\begin{subfigure}[t]{0.49\textwidth}
\centering
\begin{tikzpicture}[font=\footnotesize]
\begin{axis}[
width=0.95\textwidth,
height=1.0\textwidth,
xmin=0.5,
xmax=5.5,
ymin=0,
ymax=20,
ylabel={$\lvert \text{Iterations} \rvert$},
xtick={1, 2, 3, 4, 5},
xticklabels={1e-4, 1e-3, 1e-2, 1e-1, 1},
xlabel={$\alpha$},
axis background/.style={fill=white}
]
\addplot [mycolor1, mark=o, mark size=2.0pt, mark options={solid, mycolor1}]
table[row sep=crcr]{
1	18\\
2	10\\
3	5\\
4	3\\
5	2\\
};
\addplot [mycolor2, mark=*, mark size=1.0pt, mark options={solid, mycolor2}]
table[row sep=crcr]{
1	18\\
2	10\\
3	5\\
4	3\\
5	2\\
};
\end{axis}
\end{tikzpicture}
\caption{Number of iterations.} \label{figure:OC_2annulus_3}
\end{subfigure}
\hfill
\begin{subfigure}[t]{0.49\textwidth}
\centering
\begin{tikzpicture}[font=\footnotesize]
\begin{axis}[
width=0.95\textwidth,
height=1.0\textwidth,
xmin=0.5,
xmax=5.5,
ymode=log,
ymin=700,
ymax=4500,
yminorticks=true,
ylabel={$s$},
xtick={1, 2, 3, 4, 5},
xticklabels={1e-4, 1e-3, 1e-2, 1e-1, 1},
xlabel={$\alpha$},
axis background/.style={fill=white}
]
\addplot [mycolor1, mark=o, mark size=2.0pt, mark options={solid, mycolor1}]
table[row sep=crcr]{
1	2575.962861\\
2	1779.872548\\
3	1168.80908\\
4	936.843578\\
5	812.892827\\
};
\addplot [mycolor2, mark=*, mark size=1.0pt, mark options={solid, mycolor2}]
table[row sep=crcr]{
1	3.972042834000000e+03\\
2	2.908906343000000e+03\\
3	1.885239213000000e+03\\
4	1.384733011000000e+03\\
5	1.195165705000000e+03\\
};
\end{axis}
\end{tikzpicture}
\caption{Total time in seconds $s$.} \label{figure:OC_2annulus_4}
\end{subfigure}
\par\smallskip
\hspace{3em}
\centering
\begin{subfigure}[t]{0.9\textwidth}
\centering
\begin{tikzpicture}[font=\footnotesize]
\begin{axis}[
hide axis,
xmin=10,
xmax=50,
ymin=0,
ymax=0.4,
legend columns=2, 
legend style={draw=white!15!black,legend cell align=left, /tikz/every even column/.append style={column sep=0.4cm}}
]
\addlegendimage{mycolor1, mark=o, mark size=2.0pt, mark options={solid, mycolor1}}
\addlegendentry{conforming};
\addlegendimage{mycolor2, mark=*, mark size=1.0pt, mark options={solid, mycolor2}}
\addlegendentry{nonconforming};
\end{axis}
\end{tikzpicture}
\end{subfigure}
\caption{Stability of the method depending on the parameter $\alpha$ on the multi-patch geometry with $2$ annuli shown in~\ref{figure:2annulus_conf_1} for the conforming and not nonconforming case.} \label{figure:OC_2annulus_all} 
\end{figure}
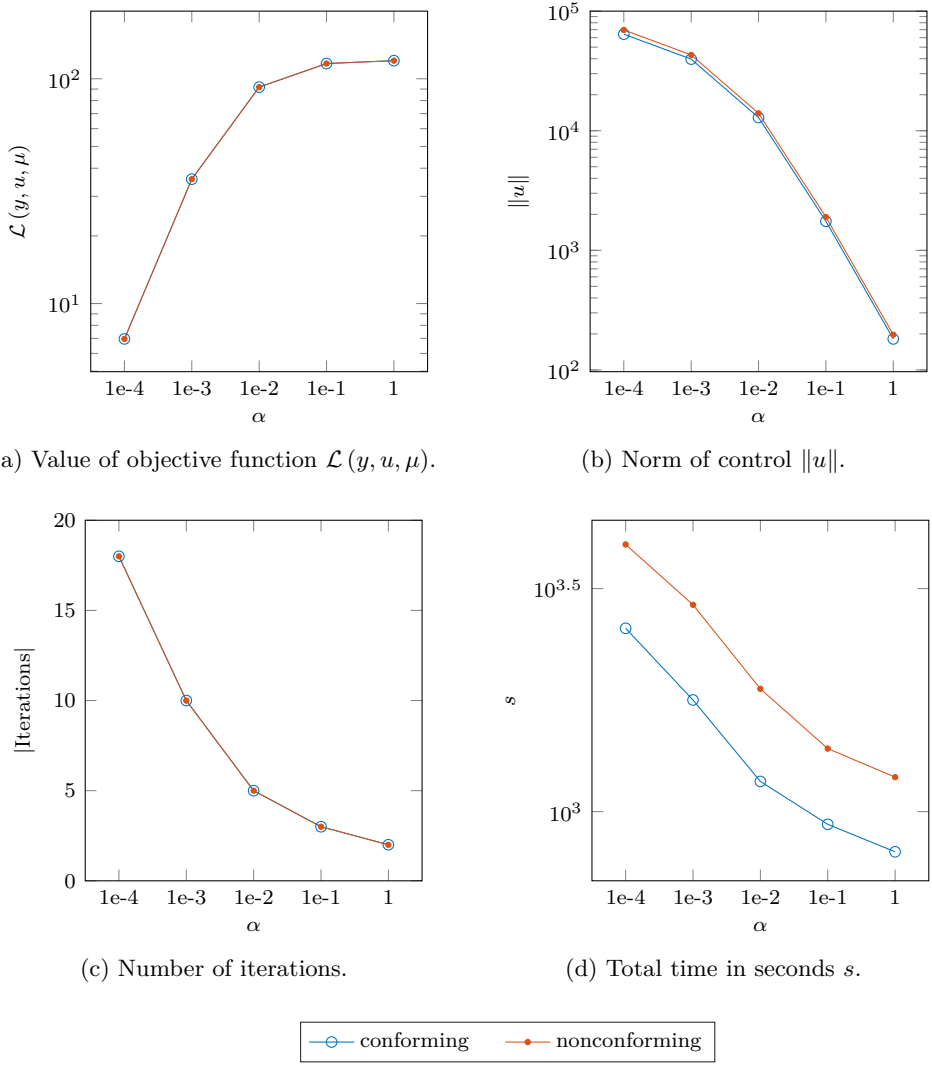

We show in Figure \ref{figure:OC_2annulus_all} the results for the multi-patch geometry shown in Figure \ref{figure:2annulus_conf_1} with desired state defined by~\eqref{eq:OC_sol} where $y_{sol} \left( x, y, z \right)$ is defined by~\eqref{eq:2annulus_sol}. We use B-splines of degree $3$ in both the conforming case and the nonconforming case. In the conforming case we have a discretization of $N^{\left( j \right)}_{DoFs} = 438.976$ DoFs for $j = 1, 2$, so that we consider a total of $8.779.520$ many DoFs for $N_t = 10$ time steps. In the nonconforming case we have the division $N^{\left( 1 \right)}_{DoFs} = 456.533$, $N^{\left( 2 \right)}_{DoFs}= 64.000$, so that we consider a total of $5.205.330$ DoFs with $N_t = 10$. We observe similar results here.

\section{Conclusion}
\label{section:conclusion}
In this paper, we transferred the TT low-rank method presented by B\"unger et al.\ \cite{BuengerDolgovStoll:2020} to the multi-patch setting using the idea of the IETI method from \cite{KleissPechsteinJuettlerTomar:2012}. The $C^0$-continuity of the approximation across the patch interfaces was ensured by defining a jump tensor which can be represented in TT format. The resulting linear system is highly structured and we showed that the solution can be approximated using TT-based solvers that rely on special \gmres\ iteration and the design of efficient but also easy to use preconditioners. We applied the resulting scheme to solve large-scale optimal control problems, where we introduced a preconditioned \gmres\ method that can deal with differently sized PDE and constraint blocks. We also equipped the method with a preconditioner that allowed for a robust performance of our scheme. We then illustrate the performance of our both methods on several multi-patch testcases.
\bibliographystyle{siam}
\bibliography{main.bbl}

\begin{thebibliography}{10}

\bibitem{Antolin}
{\sc P.~Antolin, A.~Buffa, F.~Calabr{\'o}, M.~Martinelli, and G.~Sangalli},
  {\em Efficient matrix computation for tensor-product isogeometric analysis:
  The use of sum factorization}, Comp. Method. Appl. M., 285 (2015), pp.~817 --
  828.

\bibitem{bdos-sb-2016}
{\sc P.~Benner, S.~Dolgov, A.~Onwunta, and M.~Stoll}, {\em Low-rank solvers for
  unsteady {S}tokes--{B}rinkman optimal control problem with random data},
  Comput. Method. Appl. M., 304 (2016), pp.~26--54.

\bibitem{saddlePoint}
{\sc M.~Benzi, H.~G. Golub, and J.~Liesen}, {\em Numerical solution of saddle
  point problems}, Acta Numerica, 14 (2005), pp.~1--137.

\bibitem{BuengerCode}
{\sc A.~B\"{u}nger}, {\em Low-rank tensor method for isogeometric analysis},
  2020.
\newblock
  \href{https://www.tu-chemnitz.de/mathematik/wire/codes.php}{tu-chemnitz.de/mathematik/wire/codes.php},
  Accessed: 2024-05-10.

\bibitem{BuengerDolgovStoll:2020}
{\sc A.~B\"{u}nger, S.~Dolgov, and M.~Stoll}, {\em A low-rank tensor method for
  {PDE}-constrained optimization with isogeometric analysis}, SIAM Journal on
  Scientific Computing, 42 (2020), pp.~A140--A161.

\bibitem{mlsvd}
{\sc L.~De~Lathauwer, B.~De~Moor, and J.~Vandewalle}, {\em A multilinear
  singular value decomposition}, SIAM Journal on Matrix Analysis and
  Applications, 21 (2000), pp.~1253--1278.

\bibitem{desilva-2008}
{\sc V.~de~Silva and L.-H. Lim}, {\em Tensor rank and the ill-posedness of the
  best low-rank approximation problem}, SIAM J. Matrix Anal. Appl., 30 (2008),
  pp.~1084--1127.

\bibitem{Dolgov:2013}
{\sc S.~Dolgov}, {\em {TT-GMRES}: Solution to a linear system in the structured
  tensor format}, Russian Journal of Numerical Analysis and Mathematical
  Modelling, 28 (2013).

\bibitem{ds-navier-2017}
{\sc S.~Dolgov and M.~Stoll}, {\em Low-rank solution to an optimization problem
  constrained by the {N}avier--{S}tokes equations}, SIAM J. Sci. Comput., 39
  (2017), pp.~A255--A280.

\bibitem{amen}
{\sc S.~V. Dolgov and D.~V. Savostyanov}, {\em Alternating minimal energy
  methods for linear systems in higher dimensions}, SIAM J. Sci. Comput., 36
  (2014), pp.~A2248--A2271.

\bibitem{FEM}
{\sc H.~C. Elman, D.~J. Silvester, and A.~J. Wathen}, {\em Finite elements and
  fast iterative solvers: with applications in incompressible fluid dynamics},
  Numerical Mathematics and Scientific Computation, Oxford University Press,
  Oxford, second~ed., 2014.

\bibitem{Herzog:2014}
{\sc R.~Herzog and O.~Rheinbach}, {\em {FETI-DP} methods for optimal control
  problems}, Lecture Notes in Computational Science and Engineering, 98 (2014),
  pp.~387--395.

\bibitem{CAD}
{\sc T.~Hughes, J.~Cottrell, and Y.~Bazilevs}, {\em {Isogeometric analysis:
  CAD, finite elements, NURBS, exact geometry and mesh refinement}}, Comp.
  Methods Appl. Mech. Eng., 194 (2005), pp.~4135--4195.

\bibitem{Hughes}
{\sc T.~Hughes, A.~Reali, and G.~Sangalli}, {\em Efficient quadrature for
  {NURBS}-based isogeometric analysis}, Comp. Methods Appl. Mech. Eng., 199
  (2010), pp.~301 -- 313.

\bibitem{KleissPechsteinJuettlerTomar:2012}
{\sc S.~Kleiss, C.~Pechstein, B.~Jüttler, and S.~Tomar}, {\em {IETI –
  Isogeometric Tearing and Interconnecting}}, Computer Methods in Applied
  Mechanics and Engineering, 247-248 (2012), pp.~201--215.

\bibitem{angelos2}
{\sc A.~Mantzaflaris, B.~J{\"u}ttler, B.~N. Khoromskij, and U.~Langer}, {\em
  {Matrix generation in isogeometric analysis by low rank tensor
  approximation}}, in Curves and Surfaces: 8th International Conference, Paris,
  France, June 12-18, 2014, Revised Selected Papers, Springer International
  Publishing, 2015, pp.~321--340.

\bibitem{angelos1}
\leavevmode\vrule height 2pt depth -1.6pt width 23pt, {\em Low rank tensor
  methods in {G}alerkin-based isogeometric analysis}, Comp. Methods Appl. Mech.
  Eng., 316 (2017), pp.~1062--1085.

\bibitem{Mantzaflaris_space_time}
{\sc A.~Mantzaflaris, F.~Scholz, and I.~Toulopoulos}, {\em Low-rank space-time
  decoupled isogeometric analysis for parabolic problems with varying
  coefficients}, Comp. Methods Appl. M.,  (2018 in press).

\bibitem{montardini2023lowsp}
{\sc M.~Montardini, G.~Sangalli, and M.~Tani}, {\em A low-rank isogeometric
  solver based on {Tucker} tensors}, Computer Methods in Applied Mechanics and
  Engineering, 417 (2023), p.~116472.

\bibitem{montardini2024lowrank}
{\sc M.~Montardini, G.~Sangalli, and M.~Tani}, {\em A low-rank solver for
  conforming multipatch isogeometric analysis}, 2024.

\bibitem{nguyen2015isogeometric}
{\sc V.~P. Nguyen, C.~Anitescu, S.~P. Bordas, and T.~Rabczuk}, {\em
  Isogeometric analysis: an overview and computer implementation aspects},
  Mathematics and Computers in Simulation, 117 (2015), pp.~89--116.

\bibitem{Oseledets:2011}
{\sc I.~Oseledets}, {\em Tensor-train decomposition}, SIAM J. Scientific
  Computing, 33 (2011), pp.~2295--2317.

\bibitem{osel-tt-2011}
{\sc I.~V. Oseledets}, {\em Tensor-train decomposition}, SIAM J. Sci. Comput.,
  33 (2011), pp.~2295--2317.

\bibitem{tt-toolbox}
{\sc I.~V. Oseledets, S.~Dolgov, V.~Kazeev, D.~Savostyanov, O.~Lebedeva,
  P.~Zhlobich, T.~Mach, and L.~Song}, {\em {TT-Toolbox}}, 2011.
\newblock https://github.com/oseledets/TT-Toolbox.

\bibitem{DoOs-dmrg-solve-2011}
{\sc I.~V. Oseledets and S.~V. Dolgov}, {\em Solution of linear systems and
  matrix inversion in the {TT}-format}, SIAM J. Sci. Comput., 34 (2012),
  pp.~A2718--A2739.

\bibitem{pan2020fast}
{\sc M.~Pan, B.~J{\"u}ttler, and A.~Giust}, {\em Fast formation of isogeometric
  {Galerkin} matrices via integration by interpolation and look-up}, Computer
  Methods in Applied Mechanics and Engineering, 366 (2020), p.~113005.

\bibitem{pan2021efficient}
{\sc M.~Pan, B.~J{\"u}ttler, and A.~Mantzaflaris}, {\em Efficient matrix
  assembly in isogeometric analysis with hierarchical {B-splines}}, Journal of
  Computational and Applied Mathematics, 390 (2021), p.~113278.

\bibitem{stoll2}
{\sc J.~W. Pearson, M.~Stoll, and A.~J. Wathen}, {\em Regularization-robust
  preconditioners for time-dependent {PDE}-constrained optimization problems},
  SIAM J. Matrix Anal. Appl., 33 (2012), pp.~1126--1152.

\bibitem{NURBS}
{\sc L.~Piegl}, {\em On {NURBS}, a survey}, IEEE Computer Graphics and
  Applications, 11 (1991), pp.~55--71.

\bibitem{piegl}
{\sc L.~Piegl and W.~Tiller}, {\em The NURBS Book}, Monographs in Visual
  Communication, Springer Berlin Heidelberg, 1996.

\bibitem{CPD}
{\sc M.~Sorensen, D.~Lathauwer, P.~Comon, S.~Icart, and L.~Deneire}, {\em
  Canonical polyadic decomposition with a columnwise orthonormal factor
  matrix}, SIAM Journal on Matrix Analysis and Applications, 33 (2012),
  pp.~1190--1213.

\bibitem{Spink:NURBS}
{\sc D.~Spink}, {\em Nurbs toolbox}.

\bibitem{stoll1}
{\sc M.~Stoll and T.~Breiten}, {\em A low-rank in time approach to
  {PDE}-constrained optimization}, SIAM J. Sci. Comput., 37 (2015),
  pp.~B1--B29.

\bibitem{strang}
{\sc G.~Strang and G.~Fix}, {\em {An Analysis of the Finite Element Method}},
  Wellesley-Cambridge Press, 2008.

\bibitem{geoPDEs}
{\sc R.~V{\'a}zquez}, {\em A new design for the implementation of isogeometric
  analysis in {O}ctave and {M}atlab: {G}eo{PDE}s 3.0}, Computers \& Mathematics
  with Applications, 72 (2016), pp.~523--554.

\end{thebibliography}
\end{document}